\documentclass[a4paper]{article}

\usepackage[utf8]{inputenc}
\usepackage[english]{babel}
\usepackage{latexsym,amssymb,amsthm,xcolor}
\usepackage{mathtools}
\usepackage{graphicx}
\usepackage{bbm}
\usepackage{framed}
\usepackage{fancyhdr}
\usepackage{url}

\usepackage{pdflscape}
\usepackage{hyperref}

\newcommand{\Levy}{L{\'e}vy}

\newcommand{\ind}[1]{\mathbbm{1}_{\{#1\}}}

\newcommand{\wt}{\widetilde}
\newcommand*\diff{\mathop{}\!\mathrm{d}}
\newcommand{\wh}{\widehat}
\newcommand{\noi}{\noindent}

\newtheorem{proposition}{Proposition}[section]
\newtheorem{corollary}[proposition]{Corollary}

\newtheorem{theorem}{Theorem}[section]

\theoremstyle{definition}
\newtheorem{definition}{Definition}[section]

\newtheorem{xx}{\bf xxx}

\theoremstyle{newremark}
\newtheorem{remarkx}{Remark}[section]

\newenvironment{remark}
  {\pushQED{\qed}\remarkx}
  {\popQED\endremarkx}

\usepackage{enumerate}

\renewcommand{\epsilon}{\varepsilon}
\DeclareMathOperator{\supp}{supp}

\allowdisplaybreaks[4]

\newcommand{\mcB}{\mathcal{B}}

\newcommand{\mcF}{\mathcal{F}}
\newcommand{\mcG}{\mathcal{G}}

\newcommand{\mcL}{\mathcal{L}}
\newcommand{\mcM}{\mathcal{M}}
\newcommand{\mcP}{\mathcal{P}}
\newcommand{\mcU}{\mathcal{U}}

\newcommand{\mfC}{\mathfrak{C}}

\newcommand{\mfK}{\mathfrak{K}}

\newcommand{\mfM}{\mathfrak{M}}

\newcommand{\mfu}{\mathfrak{u}}
\newcommand{\mfU}{\mathfrak{U}}
\newcommand{\mfv}{\mathfrak{v}}

\newcommand{\D}{\mathbbm{D}}
\newcommand{\E}{\mathbbm{E}}
\newcommand{\G}{\mathbbm{G}}
\newcommand{\I}{\mathbbm{I}}
\newcommand{\V}{\mathbbm{V}}

\newcommand{\N}{\mathbbm{N}}
\newcommand{\R}{\mathbbm{R}}
\newcommand{\U}{\mathbbm{U}}
\newcommand{\Z}{\mathbbm{Z}}
\newcommand{\M}{\mathbbm{M}}

\newcommand{\bbS}{\mathbb{S}}

\newcommand{\uur}{\underline{\underline{r}}}

\begin{document}

\title{Stochastic evolution of genealogies of spatial populations:\\
state description, characterization of dynamics and properties}
\author{A. Depperschmidt $^1$, A. Greven$^2$}
\date{\today}
\maketitle

\begin{abstract}
  We survey results on the description of stochastically evolving
  genealogies of populations and marked genealogies of multitype
  populations or spatial populations via tree-valued Markov processes
  on (marked) ultrametric measure spaces. In particular we explain the
  choice of state spaces and their topologies, describe the dynamics
  of genealogical Fleming-Viot and branching models by well-posed
  martingale problems, and formulate the typical results on the
  longtime behavior. Furthermore we discuss the basic techniques of
  proofs and sketch as two key tools of analysis the different forms
  of duality and the Girsanov transformation.
\end{abstract}
\bigskip

\tableofcontents

\bigskip
\footnoterule
\noi
\hspace*{0.3cm}\\
{\footnotesize $^{1)}$ Department Mathematik, Universit\"at
  Erlangen-N\"urnberg, Cauerstr.~11, D-91058 Erlangen, Germany,
  \\depperschmidt@math.fau.de}\\
{\footnotesize $^{2)}$ Department Mathematik, Universit\"at
  Erlangen-N\"urnberg, Cauerstr.~11,
  D-91058 Erlangen, Germany, \\greven@math.fau.de}\\

\newpage

\newcounter{secnum}
\setcounter{secnum}{\value{section}}
\setcounter{section}{0}
\setcounter{secnumdepth}{3}
\setcounter{equation}{0}
\renewcommand{\theequation}{\text{\arabic{secnum}.\arabic{equation}}}
\numberwithin{equation}{section}

\section{Background}\label{s.back}
The aim of this contribution is to present the work of the authors and
co-workers Peter Pfaffelhuber, Rongfeng Sun and Anita Winter, on
\emph{evolving genealogies} in \emph{spatial} population models in a
systematic way and to explain its role in the context of the field of
tree-valued Markov processes which is a field with many facets which
is developing in many directions and with many applications from
population genetics to computer science. For reasons of consistency in
this survey we deviate from the original notation at various points.

In the 1970s the study of spatial population models such as voter
model, branching random walk, contact process began and the focus was
on the understanding of the longtime behavior of these systems. In
particular the occupation measures were studied, that is, the
configurations of numbers of particles/types per site was viewed as a
measure on geographic and/or types space. This process of
\emph{occupation measures} is a Markov process and can be studied as
\emph{measure valued processes} on infinite geographic spaces; see
e.g.\ \cite{AthreyaNey1972,Lig85,Du88,D93}. Nevertheless inspired by
the graphical construction often used, the results were interpreted in
terms of \emph{genealogical} information on individuals.

In this context tracing of \emph{ancestral paths} of an individual
from the current population backwards in time to its ``ancestor''
played an increasing role and was formalized by Dawson and Perkins in
\cite{DP91} via the so called \emph{historical process}, i.e., a
measure-valued process with measures on paths which were the ancestral
paths of individuals currently alive. The approach worked well for
branching random walks and the Dawson-Watanabe super processes. In the
latter case these paths encoded indeed the genealogical distance
between individuals using the information on location or type.
However, the genealogies were \emph{not} explicitly part of the state.

The techniques of treatment of historical processes used infinite
divisibility. Conceptually the historical processes of the ancestral
paths are also of importance in other population models as in the
voter model, contact process or Moran and Fleming-Viot models, where
results could be obtained using the technique of the
\emph{hierarchical mean-field limit}, see \cite{DGV95} or latter using
ideas of the \emph{lookdown construction} of Donnelly and Kurtz
\cite{DK99a} and \emph{particle system techniques} of Liggett and
Spitzer \cite{LS81} in \cite{GLW05}. In the Fleming-Viot model on
$\Z^1$ the construction of the Brownian web \cite{Arr81,NRS05,FINR04}
can be used to obtain the processes of ancestral path, \cite{GSW}. The
historical process can be combined with the coding of genealogical
information as genealogical distances in quite some generality, even
in case of branching processes and super random walks getting again
(generalized) branching processes, provided that the genealogical
information is coded suitably \cite{ggr_GeneralBranching}. We will
discuss this coding later on.

The genealogical point of view gained more attention in the analysis
of infinite particle systems which exhibit a \emph{regime of
  clustering} in low dimensions. For example in the case of voter
model there are growing clumps of a single opinion as time diverges,
an opinion which originates from a single individual with that
opinion. The formation of such clusters appears in Fleming-Viot models
and in branching models as clumps of survivors in the sea of
unpopulated sites. In many papers the growth in time and space of such
``mono-ancestor'' clumps was studied in detail starting from Cox and
Griffeath \cite{CG86} with \cite{DG93,FG94,Kle96,Kle97,Win02}. In this
series of works it became more and more apparent that the
\emph{genealogy} of the population and their evolution is behind all
the effects.

At this point Aldous introduced in \cite{Ald1991a,Ald1991,Ald1993} the
\emph{continuum random tree} which made the genealogical distances in
a branching population explicit in the state description. The approach
focused on all individuals ever alive and was working with
\emph{metric trees}. In this context the description of the genealogy
based on an embedding of the genealogical trees in \emph{Brownian
  excursions}, such as in \cite{NP89,LG89} for example, gave an
impulse to the ideas in the theory on continuum random trees. More
systematically the concept of metric spaces, more precisely
ultrametric spaces, was applied by Evans in \cite{Ev00} to study the
genealogy associated with the entrance law from an infinite population
of the Kingman coalescent describing the individuals currently alive
with their genealogy in a Fleming-Viot population. This was taken up
by Le~Gall in \cite{LeGall93,LG99} in the work on the \emph{Brownian
  snake}. Here the location marks on the path of descent are
constructed as well. However, here also the population alive at some
time was modeled as marked labeled tree.

The ideas of tracing ancestries and genealogical relations was a key
point in the work of Donnelly and Kurtz \cite{DK99b,DK99a}. In their
\emph{lookdown construction} the evolution of the population in, for
example a Fleming-Viot model, is represented by the evolution of a
countably infinite ordered population. The state of the original
population is obtained through the state of the ordered population and
the latter is constructed by breaking the symmetry of the dynamics.

Also taking up the work of Donnelly and Kurtz \cite{DK99a,DK99b},
another method to code genealogical information of an evolving
population was invented by Bertoin and LeGall in
\cite{BLG00,BLG03,BLG05}, generalized by Dawson and Li in \cite{DL12}
and studied further in a number of concrete cases in
\cite{Fou2012,L14,Gufler2018}. This method leads to so called
\emph{``flows of bridges''} which might be seen as a specific
technique to construct an \emph{explicit representation} of a
\emph{labeled} tree which models the evolution of the genealogy of an
evolving population.

In the program ``Genealogies of interacting particle systems'' in
Singapore 2017 a learning session ``Tree-valued Markov processes''
(held by the authors of this contribution) was devoted to flows of
bridges. Some other learning sessions where devoted to Aldous' work on
continuum random trees and variants of the lookdown construction of
Donnelly and Kurtz. The approach via flows of bridges works extremely
nicely in the case of non-spatial models whereas introducing space is
a very difficult matter. We don't have the space to spell out the
details in this survey.

In 2009 Greven, Pfaffelhuber and Winter started to describe the
\emph{evolution of genealogies} of the individuals alive at the
current time $t$ of populations in the large populations limit
\cite{GPW09,GPWmp13}. This was continued in a series of papers
\cite{DGP11,DGP12,GSW,infdiv,ggr_tvF14}. This was motivated first of
all by trying to understand some effects in the longtime behavior of
spatially interacting systems and the growth in space of monotype
clusters in such systems as time diverges, a phenomenon discussed
above. Further motivation came from the interplay between the
mechanisms of migration and selection in Fleming-Viot models. As these
models are infinite population models, where also the representation
by lookdoown constructions gets a bit intransparent, we saw the
necessity to work systematically in a framework with \emph{martingale
  problems} on state spaces encoding the genealogies to an extent
needed to obtain the observed effects. To understand the genealogy of
voter models exhibiting the monotype clusters in the scaling limit on
spatial scale $\sqrt{t}$ starting with \cite{Arr81}, later
\cite{FINR04,SSS15}, went through a long struggle. We think that this
model now can be treated best in the framework of processes with
values in ultrametric measure spaces using the tools of the theory on
the Brownian web \cite{Arr81,FINR04,NRS05,SSS15}; see \cite{GSW}.

The description of the genealogy of a population and its dynamics is
by means of \emph{well-posed martingale problems} with values in the
\emph{equivalence classes of ultrametric measure spaces} respectively
\emph{marked}-ultrametric measures spaces in the case of multitype
and/or spatial populations. Using methods of stochastic calculus
\emph{path properties}, longtime behavior and properties of
\emph{equilibria} could be proved for population models driven by
\emph{Fleming-Viot resampling} or \emph{Feller branching}, i.e.\ in
the case of two large classes of models which arise as limits of
individual based models.

In the study of the evolution of genealogies of \emph{branching}
processes the idea of infinite divisibility plays also on the level of
genealogies an important role, as it does for branching processes on
$\R$. On $\R$ infinite divisibility was based on the concept of
convolution and hence it makes use of the \emph{semigroup} structure
of $\R$, which is known to be sufficient to derive \Levy{}-Khintchine
formulas. Thus, on the state space for the random genealogies one
needs to have similar concepts and therefore we need here a binary
operation with suitable algebraic properties, to be able to obtain the
\emph{\Levy{}-Khintchine} and \emph{Cox point process representations}
of the random genealogies in branching populations. Algebraic and
topological methods play an important role since for the genealogies
one can define \emph{semigroup} operations, see here Evans and
Molchanov \cite{EM17} respectively collections of such different
operations introduced by Gl\"ode, Greven and Rippl in \cite{infdiv},
\cite{ggr_GeneralBranching}, with the \emph{concatenation} as binary
operation of \emph{truncated} trees, truncation generating a
collection indexed in the truncation level, leading to a consistent
structure of connected semigroups.

We will explain in this survey how states are modeled in the
approaches mentioned in the last three paragraphs, introduce the
martingale problems for the most basic processes and describe their
properties. We consider \emph{tree-valued Fleming-Viot} and
\emph{tree-valued Feller diffusions} and later the corresponding
\emph{spatial} versions. More precisely, we consider the genealogical
processes associated with the following classical diffusion models for
populations which are the many individuals/small mass limits of
individual based models:
\begin{itemize}
\item the \emph{Fleming-Viot diffusion}, cf.\ \cite{GPWmp13},
\item the \emph{branching Feller diffusion}, cf.\ \cite{ggr_tvF14},
\item the \emph{selfcatalytic branching diffusion}, cf.\ \cite{Gl12},
\item \emph{Fleming-Viot diffusion with selection and
    mutation}, cf.\ \cite{DGP12},
\item \emph{interacting Fleming-Viot diffusions} on countable
  geographic spaces, cf.\ \cite{GSW},
\item \emph{interacting branching Feller diffusions}, cf.\
  \cite{ggr_tvF14}, and \emph{interacting logistic Feller-branching}
  diffusions, cf.\ \cite{gmuk}.
\item \emph{interacting Fleming-Viot diffusions} on the continuum,
  i.e.\ $\R$, cf.\ \cite{GSW}.
\end{itemize}
The work in this directions continues; see
\cite{Loehr13,KL15,LVW15,Gufler2018,ALW16}. Furthermore, the approach
can be extended to genealogies of individuals ever alive up to the
current time horizon. This leads to the study of the evolution of
genealogies of \emph{fossils}, cf.\ \cite{GSWfoss},
\cite{ggr_GeneralBranching}.

We choose here to outline two directions of our current research,
namely first \emph{logistic branching} processes to model populations
in competition for limited resources and the effects on genealogies
and second Fleming-Viot models with \emph{recombination}. We will
explain therefore at the end in Section~\ref{s.outlook} how
\emph{infinitely divisible} random genealogies can be studied and how
the mechanism of \emph{recombination} can be incorporated to give some
perspectives.

We would like to mention that the field develops also in various other
directions with further types of evolutions requiring new points of
view. Recently Kliem and Winter \cite{KW17} and Athreya, L\"ohr,
Winter \cite{ALW16,ALW17} have studied different types of evolution
than in the list above which however we cannot discuss here in more
detail.

\section{State spaces}
\label{s.statspac}

In the discussion of state spaces we start with the simplest case,
pass then to refinements and finally come to some topological issues
which will be needed later on.

\subsection{Ultrametric probability measure spaces and random
  genealogies}
\label{ss.ultra}

The key problem in passing from individual based models and their
genealogies to the ``continuum'' limit is the need to find a suitable
Polish space which can serve as the state space for the evolution of
\emph{random genealogies}. We next describe the solution to that
problem.

\paragraph{The state space $\U_1$}
The first objective is to find a description of the essential features
of the \emph{genealogy} of the population of currently alive
individuals of a population evolving in time and then to find a
\emph{Polish space} in which these states can be embedded. This will
allow to apply the standard methods of the theory of stochastic
processes, such as Markov process theory and martingale problems.

The key concept of the genealogical information is that of an
\emph{ancestor}. In particular the most recent common ancestor
(\emph{MRCA}) of two individuals allows to define the
\emph{genealogical tree distance} as twice the time back to the MRCA.
Hence, our approach is to describe the genealogy by a pair $(\wt U,r)$
where the set $\wt U$ labels the individuals of the population
currently alive and $r$ is the pseudo-ultrametric encoding the
genealogical distance. Note that an individual having offspring
generates individuals whose distance is zero at the time of their
birth. Therefore it is convenient to declare them as equivalent points
in $(\wt U,r)$, obtaining an ultrametric measure space $(U,r)$ and
encode the information, having here \emph{several individuals}, by
introducing a Borel probability measure on $(U,r)$ which is called the
\emph{sampling measure}. This measure allows to \emph{sample} typical
individuals from the population and is also useful in infinite
populations to obtain \emph{observable data} from the genealogy by
drawing finite samples. Hence we obtain the \emph{ultrametric
  probability measure space}
\begin{align}
  \label{e241}
  (U,r,\mu)
\end{align}
associated with the genealogy of the population of individuals
currently alive. We do not want to keep the individual names as part
of the description and therefore pass to \emph{equivalence classes} of
ultrametric probability measure spaces, which we denote by
\begin{align}
  \label{e246}
  \mcU = [(U,r,\mu)],
\end{align}
where two spaces $(U,r,\mu)$ and $(U',r',\mu')$ are said to be
\emph{equivalent} if there exists a bijective map
\begin{align}
  \label{e250}
  \varphi: \supp (\mu) \to \supp (\mu'),
\end{align}
with
\begin{align}
  \label{e253}
  r'\bigl(\varphi(x_1),\varphi(x_2)\bigr)=r(x_1,x_2), \quad \mu
  \text{ - a.s., \; and } \quad \varphi_\ast \mu = \mu'.
\end{align}
Here $\varphi_\ast \mu$ denotes the image measure of $\mu$ under the
map $\varphi$. Note that for finite genealogical trees with sampling
measure this equivalence relation means that any measure preserving
renaming of the vertices gives a population with the same genealogy.
We denote the set of all equivalence classes (to avoid set theoretic
paradoxes, we require $U \subseteq \R$ here, for example) by
\begin{align}
  \label{e257}
  \U_1.
\end{align}
This is the \emph{state space} for the process of genealogies of the
population alive at a given time $t$.

In order to obtain a decent state space for stochastic processes, we
have to equip $\U_1$ with a \emph{topology} $\mcG$ so that the
topological space $(\U_1,\mcG)$ is \emph{Polish}, and is therefore
suitable to accommodate limits of sequences of finite trees. All
populations we consider here arise as \emph{limits} of \emph{finite
  individual based} models.

We choose here the so called \emph{Gromov weak topology}. The idea to
define convergence, and thus $\mcG$, is as follows. Suppose we
consider for all $\mcU \in \U_1$ a sequence of \emph{finite sampled
  trees}. Then we define the sequence $(\mcU_n)_{n \in \N}$ as
convergent to $\mcU_\infty$, if $\mu^{\otimes \N}$ - a.s. all sampled
finite trees converge in the sense of finite graphs with a metric to
the corresponding object in $\mcU_\infty$ in law. In other words the
$(n \times n)$-matrix of distances induces a family of distributions
$(\nu_m^{(n)})_{m \in \N}$ on $(\R_+)^{\binom{n}{2}}$ which converges
weakly to $\nu_\infty^{(n)}$:
\begin{align}
  \label{e267}
  \nu_m^{(n)} \Longrightarrow \nu_\infty^{(n)}, \; \text{for all $n \in
  \{2,3,\dots\}$ as $m\to\infty$}.
\end{align}
This convergence defines a topology, the \emph{Gromov weak topology},
which has the following property.
\begin{theorem}[Polish state space \cite{GPW09}]
  \label{th.271} \leavevmode\\
  The space $\U_1$ equipped with the Gromov-weak-topology is Polish.
\end{theorem}

\begin{remark}[Metric on $\U_1$]
  \label{r.278}
  The metric can be introduced similarly to the Gromov-Hausdorff
  metric using the Prohorov metric of probability measures (instead of
  Hausdorff metric). More precisely, to compare $[(U,r,\mu)]$ and
  $[(U',r',\mu')]$ we consider \emph{isometric embeddings} (of
  representatives) into a third metric space, say $(Z,r_Z)$
  \begin{align}
    \label{281}
    U \hookrightarrow Z, \; U' \hookrightarrow Z
  \end{align}
  giving image probability measures
  $\mu_U^\ast, \mu_{U'}^\ast \in \mcP((Z,r_Z))$ and then consider the
  Prohorov metric $d_{\mathrm{P}} (\mu_U^\ast, \mu_{U'}^\ast)$. Then
  take the infimum over all such embeddings to define the distance.
  Since the quantity is independent of the choice of representatives
  of $[(U,r,\mu)]$ respectively $[(U',r',\mu')]$ this defines a
  \emph{metric} on $\U_1$. It can be shown that this metric generates
  the Gromov weak topology; see \cite{GPW09}.
\end{remark}

\paragraph{$\U_1$-valued random variables}
Next we turn to random $\U_1$-valued variables, which we denote by
\begin{align}
  \label{e277}
  \mfU.
\end{align}
The genealogy of a stochastically evolving population generates a
$\U_1$-valued stochastic process, which we denote by
\begin{align}
  \label{e281}
  \mfU=(\mfU_t)_{t \geq 0}
\end{align}
slightly abusing the notation. The corresponding path spaces
\begin{align}
  \label{e285}
  D \bigl([0,\infty),\U_1 \bigr) \quad\text{and}\quad C
  \bigl([0,\infty),\U_1 \bigr)
\end{align}
are again \emph{Polish spaces} and the theory of stochastic processes,
Markov processes and martingale problems can be applied.

In order to study the $\U_1$-valued stochastic process by means of
martingale problems or studying the longtime behavior we need a
description via \emph{test functions} on $\U_1$ whose expectations are
law determining. The key idea is to draw finite samples of trees and
then study these finite samples by test functions.

We consider therefore test functions of the following form: For
$\mcU=[(U,r,\mu)]\in \U_1$ we set
\begin{align}
  \label{e291}
  \Phi^{n,\varphi} (\mathcal{U}) = \int_\mcU \varphi \left((r(x_i,x_j))_{1
      \leq i < j \leq n} \right) \mu^{\otimes n} (\diff x_1,
  \ldots,\diff x_n),
\end{align}
where $\varphi \in C_b(\R^{\binom{n}{2}}, \R)$ for $n \ge 2$, and
$\varphi$ constant function for $n=1$. Alternatively, for $n\ge 2$ we
could write
\begin{align}
  \label{e296}
  \Phi^{n,\varphi} (\mathcal{U}) = \langle \nu^{(n)}_{\mathcal{U}},
  \varphi \rangle,
\end{align}
where $\nu^{(n)}(\mathcal{U})$ is the \emph{distance matrix
  distribution} of order $n$ of $\mcU$. We refer to such test
functions as \emph{polynomials}. The set of all polynomials (together
with the constant functions) will be denoted by
\begin{align}
  \label{e304}
  \Pi.
\end{align}
This set forms an \emph{algebra} of bounded continuous functions.
\begin{theorem}[Polynomials are law and convergence
  determining]\label{th.300}
  \leavevmode
  \begin{enumerate}[(i)]
  \item The set $\Pi$ is separating on $\U_1$ and
    $\{\E[F(\mfU)], F \in \Pi\}$ is law determining.
  \item The set $\Pi$ is convergence determining on the set of
    probability measures on $\U_1$.
  \end{enumerate}
\end{theorem}
For proofs of the above results see Theorem~5 in \cite{DGP11} or
Corollary~2.8 in \cite{Loehr13}.

\medskip

Now that we have introduced the polynomials let us note that the
topology on $\U_1$, as defined in \eqref{e267}, could also be defined
equivalently in terms of polynomials (see \cite{GPW09}):
\begin{align}
  \label{e579womarks}
  \mcU_n \xrightarrow{n \to \infty} \mcU \; \text{ in } \; \U_1
  \quad \text{iff} \quad \Phi(\mcU_n) \xrightarrow{n \to \infty}
  \Phi (\mcU) \; \text{ for all } \; \Phi \in \Pi.
\end{align}
This will be used next for extending and enriching the state spaces.

\subsection{Extensions: marks, finite and infinite populations,
  fossils}\label{ss.extens}

The concepts introduced in the previous section have to be extended in
two ways. First to handle \emph{spatial} and \emph{multitype}
populations we have to pass to \emph{marked} genealogies. Second to
incorporate \emph{varying population sizes} we have to allow sampling
measures which are not necessarily probability measures and which
could even have infinite total mass, at least globally.

\paragraph{Genealogies of multitype and spatial populations}
If we have a multitype population which is spatially distributed we
have to associate with the individuals types from some \emph{type set}
$\I$ and the types may influence the reproduction mechanism.
Similarly we may have a \emph{geographic space} $\G$ in which the
individuals are located, i.e., we assign them geographic locations. We
assume that the metric spaces $(\I,r_\I)$ and $(\G,r_\G)$ are complete
and separable.

We incorporate this in our formalism by considering a measurable map
\begin{align}
  \label{e353}
  \wt \kappa : U \to \V, \quad \V= \I \times \G,
\end{align}
where $\V$ is the mark space and $\wt \kappa(u)$ is the mark of
individual $x \in U$. Here $\V$ equipped with the product metric
$r_\V$ is also complete and separable. This results in the structure
\begin{align}
  \label{e359}
  (U,r,\mu,\wt \kappa)
\end{align}
describing a \emph{$\V$-marked genealogy}. We now introduce the
Borel-measure $\nu$ on $(U \times \V,r \otimes r_\V)$:
\begin{align}
  \label{e363}
  \nu=\mu \otimes \kappa,\quad
  \kappa(x,\diff v)=\delta_{\wt \kappa(x)} (\diff v), \; x \in U.
\end{align}
Then abstracting again from the individuals names (of course not their
types or locations) we consider the \emph{equivalence classes} of
$\V$-marked ultrametric probability measures
\begin{align}
  \label{e367}
  \left[(U \times \V,r\otimes r_\V,\nu)\right],
\end{align}
where we say that $(U \times \V,r\otimes r_\V,\nu)$ and
$(U' \times \V,r'\otimes r_\V,\nu')$ are equivalent if there is a map
$\varphi: U \times \V \to U' \times \V$ such that (here $\pi_U$ and
$\pi_\V$ are the projections on $U$ respectively $\V$)
\begin{align}
  \label{e371}
  \begin{split}
    \pi_U \circ \varphi \text{ is an } &  \text{ isometrie between } \supp
    ((\pi_U)_\ast \nu) \text{ and } \supp((\pi_{U'})_\ast \nu') \\
    \varphi((x,v)) & = \left((\pi_U \varphi)(x,v),v \right) \quad \nu
    \text{ - a.s., \; and } \\
    \varphi_\ast \nu & = \nu'.
  \end{split}
\end{align}
To obtain later a \emph{complete} space it is necessary to allow here
instead of mark functions more generally \emph{mark kernels $\kappa$}.
Then we denote the set of all equivalence classes of $\V$-marked
ultrametric probability measure spaces by
\begin{align}
  \label{e383}
  \U^\V_1,
\end{align}
which will be the state space for the marked genealogies.

On this space we introduce the \emph{$\V$-marked Gromov-weak topology}
as follows. Again we use the idea that finite samples converge as
finite marked trees. Let
\begin{align}
  \label{e552}
  \Pi^\V
\end{align}
be the set of polynomials of the form
\begin{multline}
  \label{e556}
  \Phi^{n,\varphi,\chi} \bigl( [ (U \times \V,r\otimes r_\V, \nu) ]\bigr) \\
  = \int_{(U \times \V)^n} \varphi (\uur)
  \xi(\underline{v}) \nu^{\otimes n} (\diff(x_1,v_1),\dots,\diff(x_n,v_n)),
\end{multline}
where $\varphi \in C_b([0,\infty)^{\binom{n}{2}},\R)$,
$\xi \in C_b(\V^n,\R)$ and $n \in \N_0$. Using these polynomials we
can again define convergence on $\U^\V_1$ by the requirement:
\begin{align}
  \label{e579}
  \mcU_n \xrightarrow{n \to \infty} \mcU \; \text{ in } \; \U^\V_1
  \quad \text{iff} \quad \Phi(\mcU_n) \xrightarrow{n \to \infty}
  \Phi (\mcU) \; \text{ for all } \; \Phi \in \Pi^\V.
\end{align}
The corresponding topology is called the \emph{$\V$-marked}
Gromov-weak topology and we obtain a \emph{Polish space} as a state
space for the random marked genealogies. It is metrizable by a metric
which generalizes the one described in Remark~\ref{r.278}; see
\cite{DGP11}.

\begin{theorem}[Polish state space]\label{th.390}\leavevmode\\
  The space $\U^\V_1$ equipped with the $\V$-marked Gromov-weak
  topology is a Polish space.
\end{theorem}
Now we can define evolving random genealogies $(\mfU_t)_{t \geq 0}$ of
spatial multitype populations as $\U^\V_1$-valued stochastic
processes. Again for $\U^\V_1$ and $\Pi^\V$ the analogue of
Theorem~\ref{th.300} holds; see \cite{DGP11}. As a prototype we will
discuss later the evolving genealogies of the spatial Fleming-Viot
process with selection.

\paragraph{Genealogies of populations with varying population sizes}
If the \emph{population size} of a population is fluctuating we have
to incorporate this size explicitly in the state description, think of
a multitype Feller diffusion for example. This means we consider now
an element
\begin{align}
  \label{e399}
  \mcU=(\bar \mcU,\wh \mcU), \quad \bar \mcU \in \R,\, \wh \mcU \in
  \U^\I_1 \quad,
\end{align}
where $\bar \mcU$ describes the population size and $\wh \mcU$ the
$\I$-marked genealogy from above (without the spatial component). We
note that now we consider the map
\begin{align}
  \label{e404}
  \mcU= (\bar \mcU,\wh \mcU) \mapsto \left[(U \times \I,r \otimes
  r_\I,\bar \mcU \nu) \right]
\end{align}
and we obtain the equivalence class of an ultrametric measure space
where the measure on $U \times \I$ is now in
$\mcM_{\rm fin}((U \times \I), \mcB(U \times \I))$. We denote this set
of states by
\begin{align}
  \label{e632}
  \U_{\rm fin}^\I.
\end{align}

Here we see that if we introduce a topology we have to settle, how we
want to handle the elements $\{(0,\wh\mcU) : \wh\mcU \in \U^\I_1 \}$.
Without marks these are extinct populations where genealogy could at
most make sense as a limiting object, genealogy at extinction. (In
fact below we shall discuss also the role of extinct types.)

Consider first the case without marks. Then we could use on
$\U_{\rm fin}$ simply the topology based on the convergence of
polynomials (we extend the definition of polynomials from \eqref{e291}
to $\U_{\rm fin}$). In this case a sequence of elements where the
population sizes converge to zero is automatically converging to the
$0$-element and in particular all elements
$\{(0,\wh \mcU) : \wh \mcU \in \U_1 \}$ are \emph{identified} with the
$0$ element, i.e., give one object.

Alternatively we could define convergence for a sequence with total
masses converging to zero as follows: \emph{the sequence
  $(\bar \mcU_n, \wh \mcU_n)$ converges if $\bar \mcU_n$ converges in
  $\R$ and $\wh \mcU_n$ converges in $\U_1$}. In the latter case we
have the product topology of $[0,\infty) \times \U_1$ and we
distinguish elements $\{[(0,\wh \mcU)] : \wh \mcU \in \U_1\}$ rather
than identifying them all with the $0$-element. \emph{On the set where
  $\bar \mcU > 0$ both topologies agree.}

We denote the two topological spaces by (cf.\
\cite{ggr_tvF14},\cite{ggr_GeneralBranching})
\begin{align}
  \label{e603}
  \U_{\rm fin} \text{ respectively } \U^\ast_{\rm fin}.
\end{align}
In $\U^\ast$ we can study the behavior of a \emph{population at
  extinction} more precisely.

In the case with \emph{marks} a suitable description is more subtle.
We assume in the following that the geographic space $\G$ is
countable. The types could be extinct temporarily during the
evolution, or geographic locations could be vacant during some time
periods and may be recolonized again during the others. If we choose
to define convergence by requiring that polynomials of the states
converge we get (we extend the definition of polynomials from
\eqref{e556}):
\begin{align}
  \label{599}
  \U^\V_{\rm fin}.
\end{align}

However, in spatial models it is often of interest to study the
genealogy on the way to extinction. What about the analogue of above
\eqref{e399} representation? We proceed as follows.

We consider the situation where for example the state at each site is
given as a pair (population size, genealogy), i.e., an element in
$\R_+ \times \U_1$ and the complete system as a measure on the
geographic space $\V=\G$ and the genealogy as an ultrametric measure
space which is an element of
$\bigl[\bigl(\bigcup_{g \in \G} U_g,r, \sum_{g \in \G} \mu_g\bigr)
\bigr]$ so that we have pairs in
$\mcM_{\rm fin} (\V) \otimes \U_{\lvert \V\rvert}$. Then the product
topology is the generalization leading to the state space
\begin{align}
  \label{e683}
  \U^{\V,\ast},
\end{align}
by choosing the product topology which is Polish.

In this situation of varying population sizes we obtain the random
variables and stochastic processes
\begin{align}
  \label{e413}
  \mfU=(\bar \mfU,\wh \mfU),
\end{align}
where we can distinguish \emph{population sizes} and where
$\bar \mfU$ is $[0,\infty)$-valued or $\mcM(\V)$-valued where $\mcM$
denotes the set of Borel measures and \emph{genealogies} $\wh \mfU$ of
the evolving population.
\begin{theorem}[Polish state space]
  \label{c.616}\leavevmode\\
  In both constructions the spaces $\U^\V_{\rm fin}$ and
  $\U^{\V,\ast}_{\rm fin}$ are Polish.
\end{theorem}
Again polynomials play an important role via moments even though now
their expectations for stochastic processes do not exist in general or
are not law determining, however in case we discuss later they will,
see \cite{Gl12,DGP18}. We have no space to explain the issues in
further details here.

\paragraph{Globally infinite locally finite population sizes}
If we consider for example a population in an infinite geographic
space as for example $\Z^d$ or $\R^d$ then it is often necessary to
consider populations with an \emph{overall infinite population size}
which is however \emph{locally finite}, that is, at single sites in
$\Z^d$ or in a compact subsets of $\R^d$. We now modify our
description of the genealogy to account for this possibility and
distinguish a \emph{measure-valued part} $(\bar \mcU_t)_{t \geq 0}$
describing occupation measures of the locations and
$(\wh \mcU_t)_{t\ge 0}$ describing genealogies so that for example for
$\G$ countable we may have a collection of ultrametric probability
measure spaces and corresponding population sizes, which we view again
as a measure on $\G$.

One important case is where this measure is always the counting
measure (locally) fixed population as in Fleming-Viot models, while in
super random walk these sizes fluctuate and we have to work with the
full space of locally finite measures. In other words we consider
sampling measures such that the projection on $\G$ satisfies:
\begin{align}
  \label{e700}
  \nu_\G (A) < \infty \text{ for all bonded Borel sets } A \subset
  \V.
\end{align}

We denote this space by
\begin{align}
  \label{e705}
  \U^{\V,\#}.
\end{align}

The topology is defined by considering the spaces corresponding to
bounded and closed subsets of the geographic space $\G_n \subseteq \G$
and $\G_n \uparrow \G$. We say that a sequence of elements in
$\U^{\V,\#}$ converges if all restrictions to $\G_n$ converge in the
sense described before. This defines the
\emph{Gromov-weak-$\#$-topology} which is independent of the choice of
the sequence $(\G_n)_{n \in \N}$; see \cite{GSW}.

\begin{theorem}[Polish state space]\label{th.polspace}
  \leavevmode\\
  The space $\U^{\V,\#}$ equipped with the Gromov-weak-$\#$-topology is
  Polish.
\end{theorem}
This again means that we have a suitable state space also for
\emph{spatial} multitype populations with \emph{infinite overall
  mass}.

\paragraph{Genealogies including fossils}
So far we have focused on the description of the genealogy of the
population \emph{currently alive}, leading to \emph{ultra}-metric
measure spaces. Sometimes however one is interested in the evolution
in time $t$ of the genealogy of the population alive up to time $t$
which includes the individuals alive at \emph{some time $s<t$} but
which do not necessarily have descendants at time $t$. In other words
we include the \emph{fossils} generated up to the present time $t$ and
see how this object now evolves as $t$ increases.

This will then lead to equivalence classes of metric measure spaces
denoted generically by $\M$ whose projections on subspace of
\emph{leaves} gives the object $\U$ we looked at above. Here $\U$
denotes either $\U_1$ or $\U_{\rm fin}$. The set $\M$ equipped with
the Gromov-weak topology is a Polish space which contains $\U$ as a
subspace (this holds also for the $\V$-marked case), see
\cite{GPW09,DGP11}. However in order to model the states of the
genealogy in an evolving population of all the individual which have
been alive up to the present time $T$ (and for equilibrium
considerations then also forever) it is useful to work with a subspace
of $\M^{\rm root}$ of routed versions of $\M$. By routed we mean that
to the elements of the spaces a root is added as distinguished element
which has to be preserved under forming of the equivalence classes. The
root allows to pin down time $0$ in the genealogy.

If $T \geq 0$ is the present time of the population, then the basic
set of individuals is of the form
\begin{align}
  \label{e669}
  M=\{(s,u) \mid u \in U_s, \; s \in [0,T]\} \cup \{\varrho\},
\end{align}
where $\varrho=(0,\varrho^\ast)$ is the \emph{root} used as reference
point for time and $U_s$ describes the population \emph{alive at time
  $s$}. Between the individuals alive at or before time $t$ we have a
\emph{partial order} induced by the \emph{ancestor-descendant}
relation between any two elements of the set $M$, which arise as
children of children of children etc.\ of earlier living individuals.
Such a structure arises in branching, Fleming-Viot models, in the
contact process, the voter model and many more. We obtain in such a
setup a metric measure space
\begin{align}
  \label{e674}
  (M, \wt r,\wt \mu),
\end{align}
with $\wt r(i,i')$ being the genealogical distance of the individuals
$i$ and $i'$ in $M$, which is the sum of the times back to their MRCA
and points in $U_s$ have distance $2 s$ to the root. The root allows
to define the position on the time axes uniquely. Further we think of
the population alive at time $s$, for each $s$ as equipped with a
sampling measure $\wt \mu_s$ on $(U_s,r_s)$, where $r_s$ is the
genealogical distance restricted to $U_s \times U_s$. this we lift to
a measure $\wt \mu_s$ on $[0,T]\times U_s$ by extension. We note that
with $T$ denoting the present time the individuals $(i,T)$ with
$(i \in U_T)$ play a special role as population currently alive as the
one driving the further evolution. Therefore the sampling measure of
the whole structure $M$ must be able to sample here a finite sample,
whereas from the fossils, $(s,i)$ with $s < T$ we should simply sample
according to an aggregated (over time) sampling measure from
$\wt \mu_s$, $s \in [0,T]$ which are the measures allowing to sample
from $U_s$. Hence $\wt \mu$ is chosen of the form:
\begin{align}
  \label{e682}
  \wt \mu = \int_0^T \mu_s \diff s + \mu^{\rm top},
  \quad \mu_s \text{  is a measure on  } [0,T] \times U_s
\end{align}
and $s \mapsto \mu_s$ is assumed measurable as a map
$[0,T]\to \bigcup\limits_{s \in [0,T]} [0,T] \times U_s$ supported on
$\{s\} \times U_s$, $s \in [0,T]$ and $\mu^{\rm top}$ on
$\{T\} \times U_T$.

Then we consider again the \emph{equivalence classes} of such rooted
metric measure spaces w.r.t.\ measure preserving and root preserving
isometries. We obtain for each $T$ a Polish space and the union over
all $T \geq 0$ gives
\begin{align}
  \label{e765}
  \M_{(T,\varrho)} \text{ and with marks } \M_{(T,\varrho)}^\V.
\end{align}
If we also introduce marks, we introduce the topology again via
polynomials. This is being currently worked out in \cite{GSWfoss}, see
also \cite{ggr_GeneralBranching}.

This leads to stochastic processes describing the genealogy of a
population up to time $t$, where this structure evolves then with $t$:
\begin{align}
  \label{e649}
  (\bar \mfM_t, \wh \mfM_t)_{t \in [0,T]},
\end{align}
which is a $\M_{(T,\varrho)}^\V$-valued stochastic process.

This structure allows now to describe the genealogy of a population in
time, which is an enrichment of the process $(\mfU_t)_{t \geq 0}$ with
states in $\U$, which now appear in the slices $[U_s,r_s,\mu_s]$ but
include the additional information on genealogical relations between
the individuals alive at different times \emph{before} time $T$.

Including now an infinite time horizon is possible, but needs some
thought to deal with possibly infinite distances, if we try to take
the closure of
$\mathop{\bigcup}\limits_{T \geq 0} \; \M_{(T,\varrho)}$.

\section{Martingale problems for evolving genealogies}
\label{s.martprob}

\subsection{Concepts and examples}
\label{ss.conc}

The processes of evolving genealogies are defined by \emph{well-posed
  martingale problems}. First we recall the notion of martingale
problems that we use here; see \cite{EK86}.

\begin{definition}[Martingale problem]
  \leavevmode\\
  Let \label{D:01} $E$ be a Polish space and let $\mathcal B(E)$
  denote the space of bounded measurable functions on $E$. Furthermore
  let $\mathrm P_0 \in\mathcal M_1(E)$,
  $\mathcal F \subset \mathcal B(E)$ and let $L$ be a linear operator
  on ${\mathcal B}(E)$ with domain $\mathcal F$. The law $\mathrm{P}$
  of an $E$-valued stochastic process $\mathcal X=(X_t)_{t\geq 0}$ is
  called a solution of the $(\mathrm P_0,L,\mathcal F)$-martingale
  problem if $X_0$ has distribution $\mathrm P_0$, $\mathcal X$ has
  paths in the space ${\mathcal D}_E([0,\infty))$, $\mathrm{P}$-almost
  surely, and for all $F\in\mathcal F$,
  \begin{align}
  \label{13def}
    \Big(F(X_t)-\int_0^t L F(X_s) \diff s \Big)_{t\geq 0}
  \end{align}
  is a $\mathrm{P}$-martingale with respect to the canonical
  filtration. Moreover, the $(\mathrm P_0,L,\mathcal F)$-martingale
  problem is said to be \emph{well-posed} if there is a unique
  solution $\mathrm{P}$.
\end{definition}

To define a process by means of a martingale problem we need to define
an operator $L$ acting on an appropriate subset of the set of test
functions $\mathcal{F}$. We do this in the following paragraphs for
several examples of tree-valued processes by adding more and more
evolutionary mechanisms to the processes.

\paragraph{The neutral tree-valued process without types}

The most basic evolutionary mechanisms involving genealogies are
\begin{itemize}
\item \emph{aging}, i.e., growth of genealogical distances with time;
\item the \emph{splitting} of a genealogical line in two due to
  resampling (i.e.\ birth) events.
\end{itemize}
Now the idea is to describe the evolution by describing the evolution
of all finite samples drawn by the resampling measure. This means we
consider the following test functions for the martingale problem.

For $\mathcal U = [(U,r,\mu)] \in \U_1$ consider the polynomials
\begin{align}
  \label{eq:mp-PHI}
  \begin{split}
    \Phi^{n,\varphi} (\mcU) & = \int_{U^n} \varphi
    (\underline{\underline r}) \mu^{\otimes n} (\diff \underline x)
    \\\text{with} & \quad \varphi: \R^{\binom{n}{2}} \to \R \quad
    \text{continuously differentiable},
  \end{split}
\end{align}
where $\underline x = (x_1,\dots,x_n) \in U^n$ and $\underline{\underline
  r} \coloneqq (r_{k,\ell})_{1\le k < \ell \le n}
=(r(x_k,x_\ell))_{1\le k<\ell\le n}$.

Then the above mechanisms are described by the generator
\begin{align}
  \label{e342}
  \Omega = \Omega^{\rm grow} + \Omega^{\rm res}.
\end{align}

The \emph{growth operator} acting on polynomials, see
\eqref{eq:mp-PHI} is given by
\begin{align}
  \label{e347}
  \Omega^{\rm grow} \Phi^{n,\varphi} (\mcU)
  = 2 \int_{U^n} \mu^{\otimes n} (\diff \underline{x}) \sum_{1 \le k < \ell
    \le n} \frac{\partial}{\partial r_{k,\ell}}
    \varphi(\underline{\underline{r}})
\end{align}
and the \emph{resampling generator} by
\begin{align}
  \label{e351}
  \Omega^{\textnormal{res}} \Phi^{n,\varphi} (\mcU)
  = d \int_{U^n} \mu^{\otimes n} (\diff \underline{x}) \sum^n_{1\le k < \ell \le n}
  (\varphi \circ \theta_{k,\ell} - \varphi)(\uur).
\end{align}
Here, $d>0$ is the \emph{resampling rate} and the \emph{resampling
  operator} $\theta_{k,\ell}$ replaces the distances of the $\ell$-th
component by those of $k$-th. More precisely it is defined by
$\theta_{k,\ell}\big(\underline{\underline r}) =
\underline{\underline{\tilde r}}$ and
\begin{align}
  \label{pp11b}
    \tilde r_{i,j} \coloneqq
    \begin{cases}
      r_{i,j}, & \text{ if }i,j \ne \ell, \\
      r_{i\wedge k, i\vee k} , & \text{ if } j=\ell, \\
      r_{j \wedge k, j \vee k}, & \text{ if } i=\ell.
    \end{cases}
\end{align}

The $\U_1$-valued process with generator consisting of the growth and
resampling parts, the \emph{$\U_1$-valued Fleming-Viot processes}
$\mfU^{\rm FV}$, was constructed and studied in \cite{GPWmp13}. The
approach can be extended to \emph{Cannings processes} or so called
$\Lambda$-Fleming-Viot models, see \cite{gklimw}. \emph{Mutation} and
\emph{selection} have been included in \cite{DGP12} leading to
$\U_1^\I$-valued processes. \emph{Spatial} models giving
$\U^\G$-valued Fleming-Viot processes were considered in \cite{GSW}.
Both the latter are successively explained next in two paragraphs.

Similarly one can introduce on $\U_{\rm fin}$ with these above two
operators the process $\mfU^{\rm Fel}$ describing the evolution of the
genealogies of the continuum state Feller branching diffusion and the
self-catalytic branching diffusion; see
\cite{Gl12,infdiv,ggr_GeneralBranching,ggr_tvF14}.

\paragraph{Tree-valued process with mutation and migration}

Let the type of an individual be given by its genotype and its spatial
location so that we may assume that the type space is
$\V=\G \times \I$ with $\G$ being the geographical space and $\I$ the
type space. Think for instance of $\G$ as a (finite) subset of $\Z^d$
and $\I$ as a compact space.

For $\U_1^\V$-valued processes we can take test functions of the
following form: For $\mathcal U = [(U \times \V,r\otimes r_\V,\nu)]$
set
\begin{align}
  \label{eq:Fngf}
  \Phi^{n,\varphi,\psi,\eta} (\mcU) = \int_{(U \times \G \times \I)^n}
  \varphi(\uur) \psi(\underline g) \eta(\underline u) \nu^{\otimes
  n}(\diff(\underline x,\underline g,\underline u)),
\end{align}
with $\varphi$ as in \eqref{eq:mp-PHI}, $\psi: \G^n \to \R$ and $\eta:
\I^n \to \R$, $\underline g=(g_1,\dots,g_n) \in \G^n$ and
$\underline u=(u_1,\dots,u_n) \in \I^n$.

As in the basic case, due to aging the genealogical distances grow
with time and do not affect the types or locations of individuals. Thus,
the \emph{growth operator} action on $\Phi^{n,\varphi,\psi,\eta}$ is
given by
\begin{align}
  \label{347}
  \Omega^{\rm grow} \Phi^{n,\varphi,\psi,\eta} (\mcU)
  = 2 \int_{(U\times \G\times \I)^n} \nu^{\otimes n}
  (\diff(\underline x,\underline g,\underline u)) \psi(\underline g)
  \eta(\underline u)  \sum_{1 \le k < \ell
    \le n} \frac{\partial}{\partial r_{k,\ell}}
    \varphi(\underline{\underline{r}}).
\end{align}
The \emph{resampling generator} does replace distances as before but
now it also replaces the types of individuals but not their location.
It is given by
\begin{align}
  \label{351}
  \Omega^{\textnormal{res}} \Phi^{n,\varphi,\psi,\eta} (\mcU)
  = d & \int_{(U\times \G\times \I)^n} \nu^{\otimes
  n}(\diff (\underline x,\underline g,\underline u))
\psi(\underline g) \\ & \cdot \sum_{1\le
  k < \ell \le n} \ind{g_k = g_\ell} ((\varphi \otimes \eta)
  \circ \theta_{k,\ell} - \varphi \otimes \eta) (\uur,\underline u).
                        \nonumber
\end{align}
Here, generalizing the resampling operator from \eqref{351},
$\theta_{k,\ell}$ is defined by
$\theta_{k,\ell}\big(\uur,\underline u) =
(\underline{\underline{\tilde r}},\underline{\tilde u})$ with
$\underline{\underline{\tilde r}}$ as in \eqref{pp11b} and
\begin{align}
  \tilde u_i \coloneqq
  \begin{cases}
    u_i & \text{ if }i  \ne \ell, \\
    u_k & \text{ if }i  = \ell.
  \end{cases}
\end{align}

For \emph{mutation}, let $\vartheta\ge 0$ be the \emph{mutation rate}
and let $\beta(\cdot,\cdot)$ be a Markov transition kernel on
$\I \times \I$ and set
\begin{align}
  \label{eq:oper-mut1}
  \Omega^{\rm mut} \Phi^{n,\varphi,\psi,\eta}(\mcU) = \vartheta
  \int_{(U\times \G\times \I)^n} \nu^{\otimes
  n}(\diff(\underline x,\underline g,\underline u))
  \varphi(\uur) \psi(\underline g) \sum_{k=1}^n (\beta_k \eta -
  \eta)(\underline u)
\end{align}
where
\begin{align}
  \label{eq:oper-mut2}
  (\beta_k \eta) (\underline u) =
  \int \eta(\underline u_k^v) \beta(u_k,\diff v),
\end{align}
with $\underline u_k^v  = (u_1,\dots,u_{k-1},v,u_{k+1},\dots, u_n)$.

\bigskip

For \emph{migration} let $A$ be a transition kernel on $\G \times \G$
and for $\psi: \G^n \to \R$ we set
\begin{align}
  \label{eq:Mmigr}
  (M_{g_j,g'} \psi)(g_1,\dots,g_n) =
  \psi(g_1,\dots,g_{j-1},g',g_{j+1},\dots,g_n).
\end{align}
Then,
\begin{align}
  \label{eq:oper-Migr}
  \begin{split}
    \Omega^{\rm migr} \Phi^{n,\varphi,\psi,\eta}(\mcU) & =
    \int_{(U\times \G\times \I)^n} \nu^{\otimes n}(\diff(\underline
    x,\underline g,\underline u)) \varphi(\uur) \eta(\underline u) \\
    & \qquad \cdot \sum_{k=1}^n \sum_{g' \in G} A(g_j,g') (M_{g_j,g'}
    \psi - \psi) (\underline g).
  \end{split}
\end{align}

\paragraph{Adding selection}

For \emph{selection} let $\alpha \ge 0$ be the \emph{selection
  coefficient} and let $\chi: \V \to [0,1]$ be a function that we
refer to as \emph{fitness function}. Then the selection operator is
given by
\begin{align}
  \label{eq:oper-sel}
  \begin{split}
    \Omega^{\rm sel} \Phi^{n,\varphi,\psi,\eta}(\mcU) & = \alpha
    \int_{(U\times \G\times \I)^{n+1}} \nu^{\otimes (n+1)}
    (\diff(\underline x,\underline g,\underline u)) \varphi(\uur)
    \psi(\underline g) \\ & \qquad \cdot \sum_{k=1}^n \ind{g_k = g_\ell} (\eta \cdot
    \chi_k - \eta \cdot \chi_{n+1})(\underline u).
  \end{split}
\end{align}
Here $\chi_k(\underline u) = \chi(u_k)$. And the functions on $\psi$
and $g$ are extended to domains $\R^{\binom{n+1}{2}}$ and $\G^{n+1}$.
Note that it is possible to consider more general fitness functions to
model \emph{diploid selection} and \emph{kin selection}, i.e.\ the
fitness may depend on the distance of individuals involved in the
selection event. We refer to Section~2.4.\ in \cite{DGP12} for for
background on selection that we have in mind here.

\paragraph{The branching world}
Similarly to the above we can introduce $\U_{\rm fin}$-valued
branching diffusions, for example continuous state branching processes
as $\U_{\rm fin}$-valued Feller diffusion, or the spatial analog the
$\U^\G$-valued super random walk. Here one has as generator
$\Omega^{\rm res} + \Omega^{\rm grow}$ and additionally
$\Omega^{\rm migr}$ in the spatial case acting now on polynomials
arising by allowing for finite measures instead of probability
measures in \eqref{eq:Fngf} and adding the functions
$\mfu \mapsto \bar\mfu$ and $\mfu \mapsto c$, where $c$ is a constant.
As in the Fleming-Viot model we may consider drift terms, here of
interest are sub/supercritical branching instead of just the critical
case (which is the second order part) or we may consider the analog of
selection which gives logistic branching with a \emph{competition
  term} acting only on the total masses.

\subsection{Results}
\label{ss.resex}

In this section we state some general results concerning solutions of
martingale problems and their basic properties in the case that the
linear operator is given by the ``full'' generator
\begin{align}
  \label{e342full}
  \Omega = \Omega^{\rm grow} + \Omega^{\rm res} + \Omega^{\rm mut} +
  \Omega^{\rm migr} +  \Omega^{\rm sel}
\end{align}
acting on functions from $\Pi^\V$ which are differentiable with
respect to the distance part so that the action of $\Omega^{\rm grow}$
in \eqref{347} is defined. Depending on the evolutionary forces to be
included in the considered process one can of course omit some parts
of the generator. In this case the set of polynomials $\Pi^\V$ should
also be modified accordingly. The corresponding proofs of
Theorem~\ref{thm:MP-wp} and Proposition~\ref{prop:Feller-Markov} can
be found in \cite{GPWmp13}, \cite{DGP12} and \cite{GSW} and the proof
in the case of ``full'' generator can be deduced by combination of the
relevant results in the cited papers.

\begin{theorem}[Martingale problem is well-posed]
  \leavevmode \\
  For \label{thm:MP-wp} all $\mathrm{P}_0 \in \mathcal M_1 (\U_1^\V)$, the
  $(\mathrm{P}_0, \Omega, \Pi^\V)$-martingale problem is well-posed.
\end{theorem}

A solution $\mathrm{P}$ of the $(\mathrm{P}_0,\Omega, \Pi^\V)$
martingale problem from the above theorem is referred to as
\emph{tree-valued Fleming-Viot dynamics with mutation, migration and
  selection} or \emph{$\U_1^\V$-valued Fleming-Viot dynamics} with
initial distribution $\mathrm{P}_0$ and has the following basic
properties.

\begin{proposition}[Properties of $\U_1^\V$-valued Fleming-Viot
  dynamics]
  \leavevmode\\
  Let \label{prop:Feller-Markov} $\mfU=(\mfU_t)_{t\ge 0}$ be a
  $\U_1^\V$-valued Fleming-Viot dynamics with distribution
  $\mathrm{P}$. Then the following properties hold:
  \begin{enumerate}[(i)]
  \item The sample paths of $\mfU$ are in $C_{\U_1^\V}([0,\infty))$,
    $\mathrm{P}$-a.s.
  \item For all $t>0$ we have $\mfU_t \in \U_{1,{\rm comp}}^\V$,
    $\mathrm{P}$-a.s., where $\U_{1,{\rm comp}}^\V$ refers to marked
    metric measure spaces in which the metric space components $(U,r)$
    are compact.
  \item The process $\mfU$ is a Feller and a strong Markov process.
  \end{enumerate}
\end{proposition}

The meaning of assertions (i) and (iii) in the above proposition is
self explaining. Assertion (ii) states that all genealogies except
maybe the initial one can be encoded by a compact ultra-metric space.
On the level of genealogies that means that for each $\varepsilon>0$
the population at time $t>0$ can be decomposed in finitely many
families so that the time to MRCA of individuals within each of the
families is bounded by $\varepsilon$. This should remind the reader of
the ``coming down from infinity''-property of some coalescent
processes.

The same statements as above are or will be proven on $\U_{\rm fin}$
respectively $\U_{\rm fin}^\V$ with $\V=\G$ and $\G$ a countable,
i.e., mostly infinite, Abelian group for $\U_{\rm fin}$-valued Feller
diffusion and $\U_{\rm fin}^\G$-valued super random walk in \cite{Gl12},
\cite{ggr_tvF14}, \cite{ggr_GeneralBranching} and in \cite{gmuk} for
logistic spatial branching.

\begin{remark}[Existence of solutions]\label{r.1165}
  A standard method to construct a solution to the martingale problem,
  as for example the class in \eqref{e342full} is by approximations
  with individual based models. In a number of the models appearing in
  population genetics as in \eqref{e342full} the solution to such a
  martingale problem as in \eqref{e342full} can be obtained by
  graphical construction of the ancestral relations in time, see for
  example the \emph{Brownian web} in \cite{GSW} or by \emph{lookdown
    constructions} see \cite{GLW05,gklimw,Gufler2018}, or by
  \emph{flows of bridges} constructions \cite{BLG00,Fou2012}. However
  for complicated selection or recombination models this becomes
  intransparent.
\end{remark}

\section{Long time behavior}
\label{s.longbehav}

In this section we discuss first under which conditions $\U_1$ and
more generally $\U^\V$-valued processes have \emph{unique invariant
  distributions} and are \emph{ergodic}. In presence of migration we
assume that the geographical space is \emph{finite}.

Let us consider the case without marks first. To this end we briefly
recall the \emph{Kingman coalescent measure tree}. Consider Kingman's
$N$-coalescent tree defined as follows: Starting with $N$ lines
coalesce each pair of lines which are present at rate $d$, which is
the resampling rate of the forwards in time process, and continue
until all lines have merged into one. Alternatively we can say that
when $k \in \{2,\dots,N\}$ lines are present the total coalescence
rate is $d\binom{k}{2}$ and at a coalescence event a random pair is
chosen to coalesce. Let $U_N=\{1,\dots,N\}$ be the set of leaves and
let the (random) metric $r_N$ given by the tree distance of a
Kingman's $N$-coalescent tree, that is $r_N(i,j)$ is twice the
distance to the MRCA of the individuals $i$ and $j$, $1\le i,j\le N$.
Let $\mu_N$ be the uniform measure on $U_N$ and consider the metric
measure space $\mcU_N\coloneqq [(U_N,r_N,\mu_N)]$. By Theorem~4 in
\cite{GPW09} there is an $\U_1$-valued random variable $\mcU_\infty$
so that
\begin{align}
  \label{eq:Kingm-meas-tree-conv}
  \mcU_N \xRightarrow{N\to\infty} \mcU_\infty.
\end{align}
The limiting random variable was introduced in \cite{Ev00} and is
called the \emph{Kingman coalescent measure tree}.

The following result is Theorem~3 in \cite{GPWmp13} and shows in
particular that the distribution of the Kingman coalescent measure
tree is the unique invariant distribution of the neutral Fleming-Viot
dynamics without marks.

\begin{theorem}[Long time behavior of the basic process]
  \leavevmode\\
  Let \label{ThLTbasic} $\mfU=(\mfU_t)_{t\geq 0}$ be an $\U_1$-valued
  process with $\mfU_0=\mcU \in \U_1$ and let $\mcU_\infty \in \U_1$
  be the Kingman coalescent measure tree. Then we have
  \begin{align}
    \label{eq:498bp}
    \mfU_t \xRightarrow{t\to\infty} \mfU_\infty.
  \end{align}
\end{theorem}

Now let us turn to the general case. Let $\pi_\V$ be the projection
from $X \times \V$ on $\V$. Given $\mfU_t = [(U_t,r_t,\mu_t)]$,
$t\geq 0$, we define the process $\zeta \coloneqq (\zeta_t)_{t\geq 0}$
by projecting the sampling measures to the type space, i.e., we set
\begin{align} \label{addx2}
  \zeta_t \coloneqq (\pi_\V)_\ast \mu_t, \quad t\geq 0.
\end{align}
Then $\zeta$ is the corresponding \emph{measure-valued Fleming-Viot
  process} and conditions for its ergodicity are available in the
literature. For example, under neutral evolution, i.e. $\alpha=0$,
ergodicity of $\zeta$ has been shown if the Markov pure jump mutation
process on $\I$ with has a unique equilibrium distribution; see
\cite{D93}. In the case $\alpha>0$ and $\chi \neq 0$, the process
$\zeta$ is ergodic if mutation has a parent independent component. In
the case with no parent-independent component in the mutation operator
ergodicity of $\zeta$ has been shown in \cite{EthierKurtz1998} using
coupling techniques. In \cite{DG14} a set-valued dual allows to prove
quite general ergodic theorems, even in the spatial context under some
conditions on the mutation process.

Obviously, an $\U^\V$-valued process cannot be ergodic if its
projection on the type space, i.e., the corresponding measure-valued
process, is not ergodic. Also it is easy to see that ergodicity of
$\zeta$ is implied by the ergodicity of $\mfU$. The following result
shows that also the converse implication holds.

\begin{theorem}[Long time behavior of the process with marks]
  \leavevmode\\
  Let \label{T4} $\mfU=(\mfU_t)_{t\geq 0}$ be an $\U_1^\V$-valued
  process with $\mfU_0=\mcU$ and let $\zeta$ be as above. Then, there
  exists a unique invariant distribution on $\U_{1,{\rm comp}}^\V$ and
  for a random variable $\mfU_\infty$ with this distribution for all
  $\mcU \in \U_1^\V$ we have
  \begin{align}
    \label{eq:498}
    \mfU_t \xRightarrow{t\to\infty} \mfU_\infty
  \end{align}
  if and only if $\zeta$ has a unique equilibrium distribution.
\end{theorem}

Note that Theorem~\ref{T4} tells us that ergodic distribution of types
is given by the ergodic distribution of the corresponding
measure-valued Fleming-Viot process, but does say little about the
``geometry'' of the limiting genealogies. In the neutral case with
mutation a combination of Theorem~\ref{ThLTbasic} and Theorem~\ref{T4}
provides a more complete picture.

We next discuss processes on \emph{infinitive} geographic space where
new features occur. Indeed spatial populations on infinite geographic
spaces, for example countable Abelian groups like $\Z^d$ or continuum
space models like $\R^d$ show in their longtime behavior an
interesting \emph{dichotomy} between low and high dimensions, or to be
more precise, between \emph{transient} and \emph{recurrent}
symmetrized migration arising from the difference of the positions of
two tagged ``individuals''. For example the Fleming-Viot model in its
neutral form and compact type space has for $t\to \infty$ limiting
states which are monotype if the symmetrized migration is recurrent.
This is the case for nearest neighbor migration for $d=1,2$, whereas
for the transient case, i.e., $d \geq 3$, we have equilibria with
coexistence or more precisely for every intensity measure $\theta$ on
$\I$ we have a translation invariant shift ergodic equilibrium with
this intensity. What happens on the level of genealogies?

In the transient case a common ancestor may not exist locally as we
obtain in the limit $t \to \infty$ infinite distances. We introduce
therefore the transformation of the ultrametric (into another
ultrametric) given by:
\begin{align}
  \label{e1209}
  r \mapsto r'\,,\quad r'(x,x') = 1-\exp(-r(x,x')),
\end{align}
where the distance $1$ corresponds exactly to infinite distance in the
original metric $r$. We denote the transformed process by
$(\wt \mfU_t^{FV})_{t \geq 0}$.

Then we obtain the following. Consider initial distributions of the
population on a countable (infinite) Abelian group $\G$ such that the
projection on $\V = \I \times \G$ gives a translation invariant shift
ergodic law with
$\mathrm{E}\bigl[((\pi_\V)_\ast \nu)(\{g\} \times \cdot)] = \theta (\cdot)
\in \mcM_1(\I)$. Note that this means that if $\G$ is infinite then
the ``sampling'' measure is infinite as well.

\begin{theorem}[Longtime behavior of spatial $\U^\V$-valued Fleming-Viot
  process $\mfU^{FV}$]\label{TH.1217}\leavevmode
  \begin{enumerate}[(a)]
  \item If the kernel
    $\wh a(g,g')=\frac{1}{2} \left(a(g,g')+a(g',g)\right)$ is
    recurrent then
    \begin{align}
  \label{e1220}
      \mcL[\mfU_t^{FV}] \mathop{\Longrightarrow}\limits_{t \to \infty}
      \Gamma \in \mcM_1(\U^\V),
    \end{align}
    where $\Gamma$ is concentrated on mono-ancestor and mono-type
    configurations, i.e., the states have a.s.\ a finite essential
    diameter and
    $\nu(U \times \cdot)=\int_\I (\lambda \otimes \delta_u)
    \theta(\diff u)$.
  \item For every $\theta \in \mcM_1(\I)$ there exists an equilibrium
    measure $\wt \Gamma_\theta$ which is translation-invariant and
    shift-ergodic, satisfies
    $\mathrm{E}[\nu(U \times \cdot)] = \lambda \otimes \theta$,
    $\lambda$ being the Haar measure on $\G$ normed to
    $\lambda(\{0\}) =1$, such that
    \begin{align}
      \label{e1225}
      \mcL[\wt \mfU_t^{FV}] \mathop{\Longrightarrow}\limits_{t \to
      \infty} \wt \Gamma_\theta.
    \end{align}
    The essential diameter of the states attains $1$ a.s.
  \end{enumerate}
\end{theorem}

This means that in the recurrent case we have \emph{mono-ancestor} and
\emph{mono-type} populations developing in any finite spatial window
whereas in the transient case, locally we have \emph{coexistence} of
descendants of different ancestors and, if $\theta$ is not the point
measure we also have coexistence of different types. This is proved in
\cite{gklimw} for the more general case of $\Lambda$-Fleming-Viot
models, see \cite{GSW} for information on spatial tree-valued
Fleming-Viot.

Similar questions can be addressed for \emph{spatial branching models}
where convergence to equilibria can be established. Here a
transformation of the metric as in \eqref{e1209} may be needed to
obtain equilibria in spatial models where with positive probability no
common ancestor of the population exists. See \cite{ggr_tvF14},
\cite{gmuk} for more information.

\section{Two techniques of analysis for our processes}
\label{s.techn}

There are two major tools to investigate our processes for example to
prove uniqueness or existence of martingale problems, namely
\emph{duality} and \emph{Girsanov transformation}. The former is also
crucial for analyzing the longtime behavior and the latter is useful
for obtaining path properties of more complicated processes from
corresponding path properties of simpler processes. We explain the
techniques in the following in more detail.

\subsection{Duality}
\label{ss.techndual}

We explain here the general concept of \emph{duality} and give the
results for genealogy-valued processes. Then we continue by discussing
\emph{strong duality}, a stronger concept allowing representations of
the state. Finally we discuss the so called \emph{conditional
  duality}.

\subsubsection{The concept of duality}
\label{sss.condual}
Let $(X_t)_{t \geq 0}$ be a Markov process with state space
$(E,\mcB) $, where $ E $ is a Polish space and $\mcB $ the Borel
$\sigma$-algebra. The process $ X $ is the solution of a
\emph{martingale problem} with operator $ L $ acting on the test
functions $\mcF \subseteq C_b(E,\R)$, the
$(L,\mcF,\delta_x)$-martingale problem, with $x \in E$.

Then we have a \emph{dual process} $(Y_t)_{t \geq 0}$ with state space
$(E',\mcB')$ contained in a Polish space, $\mcB'$ being the Borel
$\sigma$-algebra. This process is the solution of the
$(L', \mcF',\delta_y)$-martingale problem,
$\mcF' \subseteq C_b(E', \mcB')$, $y \in E'$. Often $L'$ and $\mcF'$
are chosen with $\mcF' \subseteq C_b(E',\R)$ such that for
$f' \in \mcF'$, $Lf' \in C_b(E',\R')$ and for the solution $Y$ the
random variables $f'(Y_t)$ are integrable for all $t$.

In order to define a \emph{Feynman-Kac duality} we define a functional
\begin{align}
  \label{e514}
  \beta_t=\int_0^t B(Y_s)\diff s
\end{align}
where $B$ is a bounded continuous function on $E'$, the so called
\emph{potential}.

\begin{definition}[Duality relation, Feynman-Kac duality]\label{def.dual}
  \leavevmode
  \begin{enumerate}[(a)]
  \item We say that $X$ and $Y$ are in \emph{duality} w.r.t.\ the
    duality function $H : E \times E' \to \R$ if the following holds:
    \begin{enumerate}[(i)]
    \item $H(\cdot,\cdot) \in C_b(E \times E',\R)$.
    \item $\{H(\cdot,y) : y \in E'\}$ is separating on $E$,
    \item the following equality holds for all $t \geq 0$:
      \begin{align}
        \label{e527}
        \mathrm{E}_x[H(X_t,y)] = \mathrm{E}_y[H(x,Y_t)], \quad \forall
        (x,y) \in E \times E'.
      \end{align}
    \end{enumerate}
  \item The process $X$ and $Y$ are in \emph{Feynman-Kac duality} if
    for all $t \geq 0$, (i) and (ii) above hold and (iii) is replaced
    by:
    \begin{align}
      \label{e532}
      \mathrm{E}_x[H(X_t,y)] = \mathrm{E}_y[H(x,Y_t) \exp(\beta_t)]
    \end{align}
    and if additionally the r.h.s.\ is integrable for all $t \geq 0$.
  \end{enumerate}
\end{definition}

\noindent The duality relation can be established by showing the
following relations:
\begin{align}
  \label{e539}
  \left( LH(\cdot,y) \right)(x)=\left(L' H(x,\cdot) \right)(y) + B(y) H(x,y),
  \quad (x,y) \in E \times E',
\end{align}
\begin{align}
  \label{e542}
  LH(\cdot,y) \in \mcF,\; L' H(x,\cdot) \in \mcF',
  \quad \forall (x,y) \in E \times E'
\end{align}
and establishing the integrability condition at the end of the
definition above.

\begin{remark}[General applications]
  The duality relations have some consequences for processes $X$ and
  $Y$ from the above definition:
  \begin{itemize}
  \item If a solution for the $(L',\mcF',\delta_y)$-martingale problem
    \emph{exists}, then a solution of the
    $(L,\mcF,\delta_x)$-martingale problem is \emph{unique}. Note that
    it is sometimes possible to drop the boundedness and/or continuity
    assumptions on the duality functions and potentials and still be
    able to use duality to establish uniqueness.

  \item If $\mcL[(X_t)_{t \geq 0}]$ is tight and if
    $\mcL[(Y_t)_{t \geq 0}] \Longrightarrow \mcL[Y_\infty]$ as
    $t \to \infty$ then
    $\mcL[X_t] \mathop{\Longrightarrow}\limits_{t \to \infty} \pi$ and
    $\pi$ is an \emph{equilibrium distribution} of the Markov process
    $X$.
  \end{itemize}
  Furthermore in many cases:
  \begin{itemize}
  \item The existence of the process $Y$ as solution to its martingale
    problem and fulfillment of \eqref{e539} and \eqref{e542} allows to
    obtain the \emph{existence} of a process $X$ solving the
    $(L,\mcF,\delta_x)$-martingale problem.
  \item The duality might be helpful to deduce the Feller property of
    processes when their duals fulfil certain properties. Roughly
    speaking this is case when the dual process is a ``finite''
    particle system. For one of the many examples see for instance
    Theorem~1(c) in \cite{DGP12} and its proof on p.~2602.
  \end{itemize}
\end{remark}

Of particular importance for us will be the case of \emph{moment
  duals}, this means for population models that the function
$H(\cdot,y)$ can be interpreted as \emph{testing} a \emph{sample of
  size $n$} from the population. In the case of the Fisher-Wright
diffusion or the Feller diffusion this gives with $E=\R$ and $E'= \N$
indeed duality functions $H(x,n)=x^n$ but the appropriate version will
hold for the case of the tree-valued Fleming-Viot process or
tree-valued Feller process.

\subsubsection{Dualities and Feynman-Kac dualities for
  \texorpdfstring{$\U_{\rm fin}$}{Ufin}-valued processes}
\label{sss.dualfeyn}
We will explain now how we obtain for the evolving genealogies of the
Fleming-Viot model respectively the Feller diffusion a \emph{moment
  duality} respectively \emph{Feynman-Kac moment duality} based on an
enriched version of the \emph{Kingman coalescent}, namely the version
\emph{enriched} by a component of \emph{distance matrices} for the
individuals. Later we explain how this carries over to the spatial
case or the multitype case with selection and mutation.

We introduce now the \emph{dual process} and the \emph{duality
  function}.

\paragraph{(1)} Consider the Kingman coalescent starting with $n$
individuals. The state of this process is an element of $\bbS_n$ the
set of all partitions of the basic set $\{1,\dots,n\}$, i.e., a tuple
$(p_1,\dots,p_{|p|})$ with $p_i \subseteq \{1,\dots,n\}$,
$i \in 1,\dots,|p|$, $\mathop{\bigcup}^{|p|}_{i=1} p_i=\{1,\dots,n\}$
and $p_i \cap p_j = \varnothing$ for $i,j \in \{1,\dots,|p|\}$ with
$i \neq j$.

We denote by $\D_n$ the set of all $n \times n$ distance matrices,
i.e., $\D_n \subset \R_+^{\binom{n}{2}}$ so that for
$(d_{i,j})_{1\le i\le j\le n} \in \D_n$ we have
$d_{i,j} \leq d_{i,k} + d_{j,k}$ and $d_{i,i}=0$ for all
$i,j,k=1,\dots,n$. We abbreviate $E'_n=\bbS_n \times \D_n$. Then the
sate space of this process is
\begin{align}
  \label{e567}
  E'=\bigcup^\infty_{n=1} E'_n.
\end{align}

The \emph{dynamics} is that every pair of partition elements
independently coalesces at rate $d$, i.e.,
$(p_1,\dots,p_{\mid p \mid}) \to (p_1,\dots, \wh p_i,\dots,p_i \cup
p_j,\dots, p_{\mid p \mid})$ where $\wh p_i$ denotes deletion at the
$i$th position. The distances between elements of the basic set grow
at speed $2$ as long as they belong to different partition elements.
The resulting process is denoted by $\mfC$. For a given $n \in \N$,
$\varphi \in C_b\big(\R_+^{\binom{n}{2}}, \R\big)$ we later want to
view $(\varphi,n)$ as part of the state of the dual. In the neutral
case the element $(\varphi,n)$ will \emph{not} change. To this end, we
extend the state space $E'$ by setting
$E'_n=\bbS_n \times \D_n \times C_b\big(\R_+^{\binom{n}{2}}, \R\big)
\times \{n\}$.

\paragraph{(2)} To define the duality function for a given $n \in \N$,
$\varphi \in C_b\big(\R_+^{\binom{n}{2}}, \R\big)$ we consider the
function
\begin{align}
  \label{574}
  H \left( [U,r,\mu], (p,(\varphi,n)) \right)
  = \int_U \dots \int_U \varphi(\uur+\uur')\diff \mu^{\otimes n}, \quad
  [U,r,\mu] \in \U_1, (p,\uur') \in E'_n.
\end{align}
We define the potential $B$
\begin{align}
  \label{579}
  B \left((p,\uur'),(\varphi,n)\right)
  = d \cdot \binom{|p|}{2}.
\end{align}

Consider now the $\U_1$-valued Fleming-Viot diffusion $\mfU^{FV}$ and
the $\U_{\rm fin}$-valued Feller diffusion $\mfU^{\rm Fel}$ with
branching rate $b$, i.e., the resampling operator has coefficient $b$,
then we have the following result.
\begin{theorem}[Duality-relation]\label{th.584}
  \leavevmode
  \begin{itemize}
  \item[(a)] The processes $\mfU^{FV}$ and $\mfC$ are in duality
    w.r.t.\ $H$ from \eqref{574}.
  \item[(b)] The process $\mfU^{\rm Fel}$ and $\mfC$ are in
    Feynman-Kac duality w.r.t.\ to $H$ and potential $B$ from
    \eqref{574} and \eqref{579}, where $d=b$.
  \end{itemize}
\end{theorem}

In the case of \emph{spatial} models we consider the spatial
coalescent where the state is enriched by a map
\begin{align}
  \label{e599}
  \Omega \to \Omega^{|p|},
\end{align}
which associates with every partition element a \emph{location}, i.e.
an element in $\G$. Accordingly the test function $\varphi$ will be
replaced by $\varphi \cdot \psi$, where $\psi$ is a bounded function
on $\G^n$, i.e. $\psi \in C_b(\G^n,\R)$. In the case of infinite
geographic spaces, i.e. on $\U^{\V,\#}$, we restrict to $\psi$ with
bounded support in $\G$. The dynamics is modified by adding a new
mechanism. The partition elements perform independent random walks
with rates $a(\cdot,\cdot)$ until they coalesce and then the
individuals in the new partition element now follow the same random
walk. Coalescence occurs only for pairs of partition elements sharing
the same geographical location.
\begin{corollary}(Duality: spatial model)
  \leavevmode\\
  \label{cor.608}
  For the $\U^\G_1$-valued interacting Fleming-Viot diffusion and
  the $\U^{\G,\#}$-valued super random walk the duality theorem above
  holds with the modified ingredients.
\end{corollary}

Similarly we can consider \emph{multitype} spatial populations with
$\V= \I \times \G$ type frequency change by resampling respectively
branching, but with the principle that types are inherited. If we have
also types we extend $\psi$, i.e.,
$\psi \in C_b \left((\I \times \G)^n, \R \right)$. In this case we
have to assign the $\I$-types to the partition element by choosing a
type for a partition element for the duality relation at time $t$ by
sampling the type from the sampling measure of the initial state of
the forward process restricted to the position of the partition
elements position at time $t$. Then partition elements will be
assigned a type when we evaluate the duality function where the type
is chosen according to the projection of $\mu$ on $\I$ at the location
$g \in \G$ at the partition element at time $T$, i.e.,
$\mu_0(U \times \{\cdot \} \times \{g\})$.

Somewhat more complicated are the multitype processes in particular if
we include \emph{selection, mutation} and \emph{recombination} in the
mechanism where already the measure-valued model requires a
\emph{function-valued dual} which is more subtle to explain already on
that level, but in fact this constructions can be carried out on
$\U^\V_1$ and even in the spatial setup.

The function-valued part is needed to have a duality once we include
\emph{selection} and \emph{mutation} or even \emph{recombination} in
the model. The dual process is built upon a spatial coalescent
enriched with a distance matrix in addition to coalescence, migration
and distance growth we have a \emph{birth} mechanism which drives a
\emph{function-valued process}, i.e.\ at the ``birth'' transition a
selection operation acts on the function. Furthermore there is action
on the function via a mutation operation. Function-valued means here
that we consider $\varphi$ and $\psi$ as part of the dual state where
$\psi$ now follows a dynamic itself. See \cite{DG14} for an extensive
treatment of the measure-valued case. We cannot give an account of
this theory here. However we show in \cite{DGP12} that one can extend
this to the $\U^\I$-valued case. The key point is to take the same
driving particle system as in the measure-valued case but to
incorporate as well the function $\varphi$ of the distances into this
picture.

\subsubsection{Strong dual representation and conditional
  duality}\label{sss.strongrep}

\paragraph{Strong dual representation}
It turns out that we can extend the dual process in many of our models
to an $\U_{\rm fin}$-valued state by passing to the \emph{entrance
  law} of the coalescent starting from countably many individuals.
Then we introduce a metric $r'(\cdot,\cdot)$ on $\N$ by defining
$T_{i,j}$ as the first time when $i,j$ are in the same partition
element if the time is less than the finite time horizon $T$.
Otherwise we set $T_{i,j}$ equal to $T$. Then we set
\begin{align}
  \label{e621}
  r'(i,j) =2T_{i,j}; \quad i,j \in \N.
\end{align}
Furthermore we define on $\N \vert_{[1,n]}$ the uniform distribution
and extend this to all subsets which have a frequency. That way we
obtain a sequence of ultrametric probability measure spaces
$[(U_n,r'_n,\mu_n)]$. These families of equivalence classes form a
tight sequence under the law of the Kingman coalescent and we get a
family of consistent laws. Therefore we have a state
$[(U', r', \mu')] \in \U_1$ on the completion of $(\N,r')$; see
Section 4 in \cite{GPW09}. We denote this random state as
$\mfC_t^{\infty,T}$.

We need next the operation of connecting two trees in a certain way
which allows to use the fact that in the evolution of
$\U_{\rm fin}$-valued processes of the type $\mfU^{\rm Fel}$,
$\mfU^{\rm FV}$ is special in the sense that the \emph{initial state}
only enters by the result of the evolution being \emph{glued} onto it.
This procedure connecting an ultrametric measure $\wt \mfU$ space onto
a given one $\mfU$, denoted $\mfU \perp \wt \mfU$, which corresponds
to grafting $\wt \mfU$ onto $\mfU$, we make precise next. Assume that
$\wt \mfU$ is the state of a process evolving for time $t$ of diameter
$2t$ from $[(U,\underline{\underline{0}},\mu)]$.

To build the glued object
$[(U, r \perp \wt r,\mu)]= [(U,r,\mu)] \perp [(\wt U,\wt r, \wt \mu)]$
we proceed as follows. We represent the ultrametric measure space
$\wt \mfU$ by embedding it into a \emph{marked weighted $\R$-tree},
such that the leaves represent $\supp(\mu)$, compare \cite{GPWmp13}
with it the $\wt \mfU$-state. Then we can define the ancestors of a
leave $\iota$ at time $0$ denotes by $\operatorname{anc}_t(\iota)$ as
the point in the tree having distance $t$ to the leave. These
ancestors we then associate with an element $\beta(\iota)$ in $U$ to
obtain the new distance function. Then we want to extend the ancestral
path of the $\R$-tree associated with $\wt \mfU$ using the one of
$\mfU$, to read off the distance $r \perp \wt r$ on $\wt U$.
Let
$(r \perp \wt r)(\iota,\iota')=\wt r
(\iota,\iota')$ if this value is less than $2t$ and
define the distances $\geq 2 t$ by
$(r \perp \wt r)(\iota,\iota')=r(\beta(\operatorname{anc}_t(\iota)),
\beta(\operatorname{anc}_t(\iota')))$. How to define $\beta (\cdot)$?

We sample from $U$ according to $\mu$ points $x_1,x_2,\dots$. Then we
want to associate these elements with the $\wt U$ founding ancestors.
How to match this with the at most countable set of all
$\wt U$-ancestors? We can for any such ancestor execute for every
element here a $\mu$-draw to match.

This works if the population has a dynamic where we can embed the
$\R$-trees associated with the states in $\U_{\rm fin}$ all into each
other. This is satisfied for the Fleming-Viot dynamics; see
\cite{Grieshammer2017}.

\begin{theorem}[Strong duality]\label{th.633} \leavevmode\\
  Let $(\mfU_t^{FV})_{t \geq 0}$ be the $\U_1$-valued Fleming-Viot
  process. Then for every $T>0$ we have
  \begin{align}
    \label{e636}
    \mcL[\mfU^{FV}_T]=\mcL[\mfU_0 \perp \mfC_T^{\infty,T}].
  \end{align}
\end{theorem}
This representation allows to generate the complete state from the
dual dynamic in form of the entrance law from a countably infinite
population and the initial state. This is the proper generalization of
the \emph{duality via graphical constructions} into our context, which
has played an important role for individual based models.

A similar result holds for the $\U_{\rm fin}$-valued Feller diffusion.
However here is an additional subtlety which we discuss next since
this is only a \emph{conditional} representation of the genealogy
$\wh \mfU$ given the path of total masses.

\paragraph{Conditioned dualities}
For many models the evolution of the total masses is again a Markov
process evolving autonomously, since the evolution does not depend
explicitly on the genealogy. In models with varying population sizes
it is therefore often possible to condition on the \emph{complete}
path $(\bar \mfU_t)_{t \geq 0}$ and to then obtain for the conditioned
law of the genealogy part $(\wh \mfU_t)_{t \geq 0}$ a dual process,
for every realization of the total mass process.

This is the case for the $\U_{\rm fin}$-valued Feller diffusion (which
solves uniquely the
$(\Omega^{\rm grow} + \Omega^{\rm res},\Pi)$-martingale problem),
where one obtains as a dual process the coalescent introduced above
but where the coalescence rate $b$ is replaced by a time-dependent
one, namely
\begin{align}
  \label{e1038}
  b \cdot (\bar \mfU_t)^{-1}, \text{  resp.  } 0 \text{ if } \bar
  \mfU_t = 0,
\end{align}
which also identifies the conditioned genealogy process as a
time-inhomogeneous Fleming-Viot diffusion; see \cite{Gl12},
\cite{ggr_tvF14}. In other words besides the Feynman-Kac duality we
also have a \emph{conditioned duality}. This means we can now get the
following analogue of \eqref{e636}.
\begin{corollary}(Conditioned strong duality: Feller)\label{cor.Feller}
  \leavevmode \\
  Consider the $\U_{\rm fin}$-valued Feller diffusion $\mfU^{\rm Fel}$
  and write
  $\mfU^{\rm Fel}=(\bar \mfU^{\rm Fel}, \wh \mfU^{\rm Fel})$, then
  \begin{align}
    \label{e1321}
    \mcL \Bigl[\wh \mfU^{\rm Fel}_T \mid (\bar \mfU^{\rm Fel}_t)_{t \geq 0}
    = \bar \mfu \Bigr] = \mcL \Bigl[\wh \mfU^{\rm Fel}_0 \perp
    \mfC_T^{\infty,T}(\bar \mfu) \Bigr].
  \end{align}
\end{corollary}
This says in particular that the process $\wh \mfU$-conditioned on the
process $\bar \mfU$ is a time-inhomogeneous $\U^1$-valued Fleming-Viot
process, so that we can obtain that the path properties of this
process carry over to the Feller cases as for example that the states
are non-atomic, or that we have a mark function if we consider
multitype models, see Section~\ref{ss.girsan} for these properties.

\subsection{Girsanov transformation}
\label{ss.girsan}

The Girsanov transformation allows us to obtain certain evolutions of
marked genealogies from some basic processes, Fleming-Viot models with
\emph{selection}, mutation and resampling from those with mutation and
resampling alone being the most important example. In that case one
can obtain the model with selection via Girsanov transform from the
one with only resampling and mutation as we shall explain below.

If the geographic space in our population models is finite or if we
even have a non-spatial model, certain of our evolutions have the
property that an additional drift term in the generator as migration,
selection, mutation or sub- and supercriticality terms in branching
models have the property that they have a \emph{density} w.r.t.\ the
process without this terms on laws in path space over some finite time
interval. In that case one can also often obtain this density
\emph{explicitly} with an appropriate \emph{Girsanov formula}. This we
explain next.

We note that on \emph{infinite} geographic space this typically fails
and one can obtain the path law only as limit of a sequence of
approximate processes on suitable \emph{finite} geographic spaces
where the above reasoning then can be applied. Next we explain this
for the $\U^\I_1$-valued Fleming-Viot process with selection.

We consider as basic process the $\U^\I_1$-valued Fleming-Viot
diffusion with mutation. Then we consider the case with selection with
fitness function $\alpha \Psi$ where $\alpha \in \R^+$. We denote the
respective processes by $\mfU^{\rm FV,\rm mut,0}$ and
$\mfU^{\rm F,\rm mut,\alpha \Psi}$.

For the Girsanov transform we need some ingredients.
Based on the fitness function $\chi'$ consider
\begin{equation}
  \label{e1685}
  \Psi(\mcU)= \frac{\alpha}{\gamma} \langle \nu^\mcU, \chi'_{1,2}
  \rangle \mbox{  where  }
  \chi'_{k,\ell}(\uur,\underline{u})=\chi'(u_k,u_\ell,r_{k \wedge
    \ell}, r_{\ell,k}).
\end{equation}
Define $\mathcal{M}=(M_t)_{t \geq 0}$
where
\begin{align}
  \label{e1355}
  M_t=\Psi(\mfU_t)-\Psi(\mfU_0)-\int^t_0 \Omega_\alpha \Psi(\mfU_s)
  \diff s, \quad t \geq 0.
\end{align}
Denoting by $[\, \cdot\, ]$ the quadratic variation of a
semimartingale we have the following result. Its proof is a
combination of Theorem~2 and equation (3.38) in \cite{DGP12}.

\begin{theorem}[Absolute continuity and Girsanov transform]\label{th.1352}
  \leavevmode
  \begin{itemize}
  \item[(a)] The law
    $\mcL \bigl[(\mfU_t^{\rm FV,\rm mut,0})_{t \in [0,T]} \bigr]$ and
    $\mcL \bigl[(\mfU_t^{\rm FV,\rm mut,\alpha \Psi})_{t \in [0,T]}
    \bigr]$ are equivalent for $\alpha \geq 0$, for $T >0$.
  \item[(b)] The Radon-Nikodym derivative is given by
    \begin{align}
    \label{e1356}
      \frac{\diff Q^\alpha_T}{\diff P_T} (\mfU)= e^{M_T-\frac{1}{2}
      [\mathcal{M}]_T},
    \end{align}
    where the quadratic variation of $\mcM$ is given by the formula:
    \begin{align}
      \label{e1571}
      [\mcM]_t
      = \frac{\alpha^2}{\gamma} \int_0^t \operatorname{Var}_{\nu_s}
      [\chi'] \,\diff s,
    \end{align}
    where we write $\nu_s = \nu^{\mfU_s}$ and
    $\operatorname{Var}_{\nu_s}[\, \cdot\, ]$ denotes variance w.r.t.\
    measure $\nu_s$.
  \end{itemize}
\end{theorem}

\begin{remark}[First and second order operators]
  The key point here is that the resampling operator is a second order
  operator whereas all other operators are first order operators; cf.\
  Section~4 in \cite{DGP12}. This allows to use here some abstract
  theory on such operators due to Bakry and \'Emery; see \cite{BE85}.
\end{remark}

The Girsanov representation has two mayor applications:
\begin{itemize}
\item reducing the \emph{well-posedness} of the processes with
  selection to the one of the basic process,
\item transferring \emph{path properties} of the basic process to the
  one with selection.
\end{itemize}

For example the following path properties can be proven for the
Fleming-Viot processes with selection by using the equivalence of laws
with the neutral case, which is simpler to analyze (see \cite{DGP12},
\cite{DGP13} for the results below):
\begin{itemize}
\item Continuity of paths.
\item The asymptotics of the number of ancestors time $\varepsilon$
  back from the present time $t$ is $2/\varepsilon$ as
  $\varepsilon \downarrow 0$.
\item The non-atomicity of the states of the process for positive
  times, i.e., the fact that with probability one for $t > 0$ the
  sampling measure projected on $U$ has no atoms.
\item The Fleming-Viot process with selection has a mark function for
  $t>0$ if this holds for $t=0$, i.e., the sampling measure on
  $U \times \I$ can be written as
  $\mu_t(\diff x,\diff v)=\mu_t^U (\diff x) \delta_{\kappa_t (x)}
  (\diff v)$ for a measurable $\kappa: U \to \I$, which is known for
  the neutral process; see \cite{KL15}.
\end{itemize}

This technique is also of use trying to handle models with
recombination and selections, since it is possible to focus first on
building recombination into the neutral model and to then add
selection, see \cite{DGP18}.

\section{Compactness and tightness of marked genealogies}
\label{s.topiss}
We discuss here two points namely what are \emph{compact} sets in the
spaces $\U_1$, $\U_{\rm fin}$ and their marked versions and how to
formulate \emph{tightness} criteria for probability measures on this
spaces, respectively processes with values in these spaces.

\paragraph{(i)}
Most basic is the criterion for compactness of a set $K$ in
$\U_1$ where one needs for $\mcU=[(U,r,\mu)] \in K$:
\begin{itemize}
\item bounded distances,
\item a no dust condition.
\end{itemize}
The first requires that pairwise distances in a set are bounded, i.e.
\begin{align}
  \label{e1331}
  R(\mcU) := \underset{\mu \otimes
  \mu}{\operatorname{ess\,sup}}[r(i,i')] \leq M, \quad \forall i,i'
  \in U, \mbox{  for all  } \mcU \in K,
\end{align}
while the second guarantees that the number $N_\varepsilon (\mcU)$ of
$\varepsilon$-balls needed such that the sampling measure of the
subset of $U$ covered by the union of the balls exceeds
$(1-\varepsilon)$ is bounded for all $\mcU \in K$ and all
$\varepsilon > 0$. Note that $N_\varepsilon(\mcU)$ can be viewed as
the number of \emph{ancestors} of an \emph{$(1-\varepsilon)$-fraction}
of the population.

In $\U_{\rm fin}$ we need in addition to the above properties that the
total masses are bounded, i.e.\ for some $M < \infty$:
\begin{equation}
  \label{1775}
  M(\mcU)= \bar \mu \leq M \quad , \quad \forall \; \mcU=[U,r,\mu] \in K.
\end{equation}

If we want to extend this to $\U^\V_{\rm fin}$ the key point is that
we need for a set $K \subseteq \U^\V_{\rm fin}$ to be compact that the
projections on $\U_{\rm fin}$ and on $\mcM_{\rm fin}(\V)$ are compact,
the latter in the weak topology. (On $\U^{\V,\#}$ one applies this to the
localizations again to get the same picture.)

\paragraph{(ii)}
Next we pass to random variables with values in $\U_1$, $\U_{\rm fin}$
and their marked versions, as well as the laws of the random
variables. For \emph{tightness} of a set of
$\mfK \subseteq \mcM_1(\U)$ of laws we need that
\begin{equation}\label{e1817}
  \{ \mcL_\Gamma \left[\left( R(\mcU),N_\varepsilon(\mcU),M(\mcU)\right)
  \right], \Gamma \in \mfK \}
\end{equation}
is for all $\varepsilon > 0$ tight in $\mcM_1([0,\infty)^3)$.

\paragraph{(iii)}
The key task is now to pass to \emph{stochastic} processes and to
obtain the tightness of laws on \emph{sets of path} of stochastic
processes from the above facts. This requires to get \eqref{e1817}
``uniformly in time'' by representing the process by graphical
constructions, lookdown constructions or using special analytic
information, on the operator of the martingale problem. Here one uses
typically more about the ancestral structure and its evolution in
time. Compare for example \cite{GPWmp13,GSW,Grieshammer2017,gklimw}.

\section{Perspectives}
\label{s.outlook}

In this section we describe some current development in this area, in
particular we describe some questions where the concepts introduced so
far need to be further refined and generalized.

\paragraph*{(1) The branching world}
One issue is to study the branching world in more detail, from
classical branching to logistic branching.

\emph{(i)} The first point is to establish some important properties
of the $\U_{\rm fin}$-valued processes of the branching type, as the
$\U_{\rm fin}$-valued version of the \emph{branching property}, a
generalized version of the concept of \emph{infinite divisibility} and
a corresponding \emph{\Levy{}-Khintchine formula} and to then use this
to obtain more specifics of the longtime behavior and structure of
populations. This point is pursued in \cite{infdiv} and now in
\cite{ggr_GeneralBranching}, \cite{ggr_tvF14}.

We introduce the subsets of $\U_{\rm fin}$ denoted
\begin{align}
  \label{e1656}
  \U_{\rm fin}(h) \text{ and } \U_{\rm fin} (h)^\sqcup
\end{align}
as the elements with essential diameter $<2h$ respectively $\leq 2h$
for $h>0$.

The basic concept which is needed here is that of the \emph{$h$-tree
  tops} of $\mcU=[(U,r,\mu)]$, denoted $\mcU(h)$ for $h>0$. Namely
consider for $(U,r)$ now $(U,r \wedge 2h)$ and restrict $\mu$ to the
corresponding $\sigma$-algebra to obtain the $h$-tree top
\begin{align}
  \label{e1707}
\mcU(h) \in \U_{\rm fin}(h)^\sqcup.
\end{align}
We can then consider the \emph{$h$-concatenation as} binary operation
in $\U_{\rm fin}(h)^\sqcup$:
\begin{align}
  \label{e1664}
  [(U_1,r_1,\mu_1)] \sqcup
  [(U_2,r_2,\mu_2)]=[(U_1 \uplus U_2,r_1 \oplus r_2,\wt \mu_1 + \wt
  \mu_2)],
\end{align}
where
\begin{align}
  \label{e1668}
  r_1 \oplus r_2(i_1,i_2)=
  \begin{cases}
    r_1 (i_1,i_2) & \text{ if } i_1,i_2 \in U_1,\\
    r_2 (i_1,i_2) & \text{ if } i_1,i_2 \in U_2,\\
    2h & \text{ if } i_1 \in U_1, i_2 \in U_2.
  \end{cases}
\end{align}
This operation defines for each $h>0$ a
\begin{align}
  \label{e1727}
  \mbox{\emph{topological semigroup} } (\U_{\rm fin}(h)^\sqcup, \sqcup^h).
\end{align}
In particular we can associate with $\U_{\rm fin}$ the collection of
\emph{semigroups} $\{(\U_{\rm fin}(h)^\sqcup, \sqcup^h),h > 0\}$ and the
collection $\{T_h,h>0\}$ of \emph{truncation maps}
\begin{align}
  \label{e1678}
  T_h:\U \to \U_{\rm fin}(h)^\sqcup, \; \mcU \mapsto \mcU(h).
\end{align}
One has acting on the semigroup in \eqref{e1727} the property that
\begin{align}
  \label{e1682}
  T_{h'} (\mfu \sqcup^h \mfv)= T_{h'}(\mfu)
  \sqcup^{h'} T_{h'} (\mfv) \text{  \emph{(consistency)}},
\end{align}
for all $h' \in [0,h),h>0$, $\mfu,\mfv \in \U_{\rm fin}(h)^\sqcup$.

Therefore we associate with $\U_{\rm fin}$ the pair of collections of
\emph{topological semigroups and truncation maps}:
\begin{align}
  \label{e1688}
  \left(\{(\U_{\rm fin}(h)^\sqcup, \sqcup^h),h>0\}, \{T_h; h>0\} \right).
\end{align}
With these concepts one can generalize the \emph{branching property}
of an $\U_{\rm fin}$-valued process. We can define \emph{infinite divisibility}
for $\U_{\rm fin}$-valued random variables and prove a \emph{\Levy{}-Khintchine
  representation} in terms of \emph{concatenations of Cox point
  processes} on $\U_{\rm fin}(h)^\sqcup$ for all $h>0$.

An $\U_{\rm fin}$-valued random variable $\mfU$ is \emph{$t$-infinitely
  divisible} if for every $n \in \N$ and $h \in [0,t]$:
\begin{align}
  \label{e1694}
  \mfU(h) \overset{(d)}{=} \mfU^{(h,1)} \sqcup^h \dots
  \sqcup^h \mfU^{(h,n)}
\end{align}
for i.i.d.\ $(U^{(h,i)})_{i=1,\dots,n}$ with values in $\U_{\rm fin}(h)^\sqcup$.

Then the \emph{\Levy{}-Khintchine formula} on $\U_{\rm fin}$ reads as follows.
The formula represents the expectation of certain exponentials via the
\emph{\Levy{} measure}. We need the following ingredients.

We define the Laplace-transform of the random variable $\mfU$:
\begin{align}
  \label{e1775}
  L_\mfU (\Phi^{n,\phi})= \mathrm{E} \bigl[\exp(-\Phi^{n,\phi} (\mfU))
  \bigr],
\end{align}
for a positive polynomial $\Phi^{n,\varphi}$.
We also need the $h$-truncated polynomial
\begin{align}
  \label{e1790}
  \Phi^{n,\phi}_h=\Phi^{n,\phi_h} \quad,\quad \phi_h(\uur)=\phi(\uur
  \wedge 2h).
\end{align}

Then the \emph{\Levy{}-Khintchine formula} for an \emph{infinitely
  divisible random variable} $\mfU$ says that there exists a unique
\emph{\Levy{} measure} $\lambda_\infty \in \mcM^\#(\U \setminus \{0\})$
with $\int (\bar \mfu \wedge 1)\lambda_\infty(\diff\mfu) < \infty$ such
that for all $h>0 (h \in (0,t]$ for $t$-inf.div.):
\begin{align}
  \label{e1781}
  -\log L_\mfU(\Phi_h)=\int_{\U_{\rm fin}(h)^\sqcup \setminus \{0\}}
  (1-e^{-\Phi_h(\mfu)}) \lambda_h(\diff\mfu),
\end{align}
where
\begin{align}
  \label{e1785}
  \lambda_h(\diff\mfu)=\int_{\U \setminus \{0\}} \lambda_\infty (d\mfv)
  1(\mfv(h) \in \diff\mfu) \in \mcM^\# (\U_{\rm fin}(h)^\sqcup \setminus \{0\}).
\end{align}
This allows to write the $h$-tops of $\mfU$ as follows:
\begin{equation}\label{e1965}
  \mfU(h)=\mathop{\sqcup}\limits_{\mfu \in N^{\lambda_h}} \mfu,
\end{equation}
where $N^{\lambda_h}$ is a Poisson point process on
$\mcM^\# \backslash \{0\}$ with intensity measure $\lambda_h$.

The marginals of the $\U_{\rm fin}$-valued Feller diffusion are infinitely
divisible if we start in an infinitely divisible initial law, see so
that we can ask for the \Levy{} measure for $\mfU$ and its properties.
This is carried out in \cite{ggr_tvF14}.

Furthermore can we ask for the laws conditioned to survive forever or
the size-biased law. These questions are studied in detail in
\cite{infdiv},\cite{ggr_GeneralBranching} and \cite{ggr_tvF14}.

\emph{(ii)} Another direction here is to consider the \emph{analog of
  selection} for branching type processes where a total mass dependent
sub-/supercriticality in the branching rate occurs which reflects the
\emph{carrying capacity $K$} of the environment where a population
lives, introducing supercriticality for states where the total
population size is below $K$ and subcriticality if it is above $K$.
This is studied in \cite{gmuk}, where the $\U_{\rm fin}^\G$-valued spatial
logistic Feller branching diffusion is constructed and where the
longtime behavior is studied via the \emph{conditioned duality}. The
goal is to understand the changes in the genealogy induced by the
competition for resources.

\paragraph*{(2) Tree-valued processes with recombination}

From the point of view of classical population genetics a natural step
in generalizing the tree-valued processes with mutation migration and
selection is to also include \emph{recombination}. The resulting
process then would involve all the important evolutionary forces. To
this end we need to extend the notion of (marked) ultra-metric measure
spaces appropriately to incorporate \emph{multi-locus genealogies}.
This is work in progress, \cite{DGP18}.

To better understand some of the difficulties that arise in the case
with recombination, let us consider a processes without mutation,
selection and migration and assume that the set of loci is given by a
finite set $\mathbbm L$. (All the forces that we exclude here for the
sake of simplicity can be included and it is also possible to consider
infinitely countable sets of loci.)

In presence of recombination with a finite set of loci $\mathbbm L$
the genealogical distance of individuals along different loci might be
different. The genealogies along all loci together can be encoded by a
\emph{multi-ultra-metric space} $(U,(r^\ell)_{\ell \in \mathbbm L})$.
In the single locus case the sampling measure $\mu$ was a measure on
the Borel-$\sigma$-algebra on $U$ induced by the corresponding metric.
In the multi-locus case the set of individuals is still given by $U$
but the sampling measure should be somehow connected to the whole set
of metrics $(r^\ell)_{\ell \in \mathbbm L}$. As the domain for the
sampling measures we choose the Borel-$\sigma$-algebra on
$U^{\mathbbm L}$ generated by the product measure induced by
$(r^\ell)_{\ell \in \mathbbm L}$ and we require that the sampling
measures are concentrated on the diagonal of $U^{\mathbbm L}$.

One can then introduce a concept of equivalence of such spaces and
pass to the set of equivalence classes. Furthermore one can introduce
a generalization of polynomials, distance matrix distributions and a
generalized notion of Gromov-weak convergence. The corresponding
topology can again be generated by a suitable Gromov-Prohorov metric
generating a Polish state space for multi-loci-genealogies.

Finally, the tree-valued processes with recombination can be obtained
as solutions to martingale problems.

\section*{Acknowledgments}

We would like to thank Matthias Birkner, Rongfeng Sun and Jan Swart
for organizing the very stimulating program ``Genealogies of
interacting particle systems'' in Singapore 2017. We thank the
Institute for Mathematical Sciences at the National University of
Singapore for hospitality. Finally, we thank the referee for careful
reading of the first version of the manuscript and helpful comments
which improved the readability.


\begin{thebibliography}{GKW18}

\bibitem[Ald91a]{Ald1991a}
David Aldous.
\newblock The continuum random tree. {I}.
\newblock {\em Ann. Probab.}, 19(1):1--28, 1991.

\bibitem[Ald91b]{Ald1991}
David Aldous.
\newblock The continuum random tree. {II}. {A}n overview.
\newblock In {\em Stochastic analysis ({D}urham, 1990)}, volume 167 of {\em
  London Math. Soc. Lecture Note Ser.}, pages 23--70. Cambridge Univ. Press,
  Cambridge, 1991.

\bibitem[Ald93]{Ald1993}
David Aldous.
\newblock The continuum random tree. {III}.
\newblock {\em Ann. Probab.}, 21(1):248--289, 1993.

\bibitem[ALW16]{ALW16}
Siva Athreya, Wolfgang L\"ohr, and Anita Winter.
\newblock The gap between {G}romov-vague and {G}romov-{H}ausdorff-vague
  topology.
\newblock {\em Stochastic Process. Appl.}, 126(9):2527--2553, 2016.

\bibitem[ALW17]{ALW17}
Siva Athreya, Wolfgang L\"ohr, and Anita Winter.
\newblock State spaces of (continuum) trees: R-trees versus algebraic trees,
  2017.
\newblock In preparation. Talk by W.~L\"ohr at {W}orkshop: {G}enealogies of
  {I}nteracting {P}article {S}ystems, 07 Aug 2017, IMS Singapore.

\bibitem[AN72]{AthreyaNey1972}
Krishna~B. Athreya and Peter~E. Ney.
\newblock {\em Branching processes}.
\newblock Springer-Verlag, New York-Heidelberg, 1972.
\newblock Die Grundlehren der mathematischen Wissenschaften, Band 196.

\bibitem[Arr81]{Arr81}
Richard~Alejandro Arratia.
\newblock Coalescing {B}rownian motions and the voter model on $\mathbb{Z}$.
\newblock Unpublished partial manuscript, 1981.

\bibitem[B{\'E}85]{BE85}
Dominique Bakry and Michel {\'E}mery.
\newblock Diffusions hypercontractives.
\newblock In {\em S\'eminaire de probabilit\'es, {XIX}, 1983/84}, volume 1123
  of {\em Lecture Notes in Math.}, pages 177--206. Springer, Berlin, 1985.

\bibitem[BLG00]{BLG00}
Jean Bertoin and Jean-Fran\c{c}ois Le~Gall.
\newblock The {B}olthausen-{S}znitman coalescent and the genealogy of
  continuous-state branching processes.
\newblock {\em Probab. Theory Related Fields}, 117(2):249--266, 2000.

\bibitem[BLG03]{BLG03}
Jean Bertoin and Jean-Fran{\c{c}}ois Le~Gall.
\newblock Stochastic flows associated to coalescent processes.
\newblock {\em Probab. Theory Related Fields}, 126(2):261--288, 2003.

\bibitem[BLG05]{BLG05}
Jean Bertoin and Jean-Fran{\c{c}}ois Le~Gall.
\newblock Stochastic flows associated to coalescent processes. {II}.
  {S}tochastic differential equations.
\newblock {\em Ann. Inst. H. Poincar\'e Probab. Statist.}, 41(3):307--333,
  2005.

\bibitem[CG86]{CG86}
J.~Theodore Cox and David Griffeath.
\newblock Diffusive clustering in the two-dimensional voter model.
\newblock {\em Ann. Probab.}, 14(2):347--370, 1986.

\bibitem[Daw93]{D93}
Donald~A. Dawson.
\newblock Measure-valued {M}arkov processes.
\newblock In {\em \'{E}cole d'\'{E}t\'e de {P}robabilit\'es de {S}aint-{F}lour
  {XXI}---1991}, volume 1541 of {\em Lecture Notes in Math.}, pages 1--260.
  Springer, Berlin, 1993.

\bibitem[DG93]{DG93}
Donald~A. Dawson and Andreas Greven.
\newblock Hierarchical models of interacting diffusions: multiple time scale
  phenomena, phase transition and pattern of cluster-formation.
\newblock {\em Probab. Theory Related Fields}, 96(4):435--473, 1993.

\bibitem[DG14]{DG14}
Donald~A. Dawson and Andreas Greven.
\newblock {\em Spatial {F}leming-{V}iot models with selection and mutation},
  volume 2092 of {\em Lecture Notes in Mathematics}.
\newblock Springer, Cham, 2014.

\bibitem[DG18]{ggr_tvF14}
Andrej Depperschmidt and Andreas Greven.
\newblock Tree-valued {F}eller diffusion.
\newblock in preparation, 2018.

\bibitem[DGP11]{DGP11}
Andrej Depperschmidt, Andreas Greven, and Peter Pfaffelhuber.
\newblock Marked metric measure spaces.
\newblock {\em Electron. Commun. Probab.}, 16:174--188, 2011.

\bibitem[DGP12]{DGP12}
Andrej Depperschmidt, Andreas Greven, and Peter Pfaffelhuber.
\newblock Tree-valued {F}leming-{V}iot dynamics with mutation and selection.
\newblock {\em Ann. Appl. Probab.}, 22(6):2560--2615, 2012.

\bibitem[DGP13]{DGP13}
Andrej Depperschmidt, Andreas Greven, and Peter Pfaffelhuber.
\newblock Path-properties of the tree-valued {F}leming-{V}iot process.
\newblock {\em Electron. J. Probab.}, 18(84):1--47, 2013.

\bibitem[DGP18]{DGP18}
Andrej Depperschmidt, Andreas Greven, and Peter Pfaffelhuber.
\newblock Evolving genealogies of {F}leming-{V}iot processes with
  recombination.
\newblock in preparation, 2018.

\bibitem[DGV95]{DGV95}
Donald~A. Dawson, Andreas Greven, and Jean Vaillancourt.
\newblock Equilibria and quasiequilibria for infinite collections of
  interacting {F}leming-{V}iot processes.
\newblock {\em Trans. Amer. Math. Soc.}, 347(7):2277--2360, 1995.

\bibitem[DK99a]{DK99b}
Peter Donnelly and Thomas~G. Kurtz.
\newblock Genealogical processes for {F}leming-{V}iot models with selection and
  recombination.
\newblock {\em Ann. Appl. Probab.}, 9(4):1091--1148, 1999.

\bibitem[DK99b]{DK99a}
Peter Donnelly and Thomas~G. Kurtz.
\newblock Particle representations for measure-valued population models.
\newblock {\em Ann. Probab.}, 27(1):166--205, 1999.

\bibitem[DL12]{DL12}
Donald~A. Dawson and Zenghu Li.
\newblock Stochastic equations, flows and measure-valued processes.
\newblock {\em Ann. Probab.}, 40(2):813--857, 2012.

\bibitem[DP91]{DP91}
Donald~A. Dawson and Edwin~A. Perkins.
\newblock Historical processes.
\newblock {\em Mem. Amer. Math. Soc.}, 93(454):iv+179, 1991.

\bibitem[Dur88]{Du88}
Richard Durrett.
\newblock {\em Lecture notes on particle systems and percolation}.
\newblock The Wadsworth \& Brooks/Cole Statistics/Probability Series. Wadsworth
  \& Brooks/Cole Advanced Books \& Software, Pacific Grove, CA, 1988.

\bibitem[EK86]{EK86}
Stewart~N. Ethier and Thomas~G. Kurtz.
\newblock {\em Markov processes: Characterization and convergence}.
\newblock Wiley Series in Probability and Mathematical Statistics: Probability
  and Mathematical Statistics. John Wiley \& Sons Inc., New York, 1986.

\bibitem[EK98]{EthierKurtz1998}
Stewart~N. Ethier and Thomas~G. Kurtz.
\newblock Coupling and ergodic theorems for {F}leming-{V}iot processes.
\newblock {\em Ann. Probab.}, 26(2):533--561, 1998.

\bibitem[EM17]{EM17}
Steven~N. Evans and Ilya Molchanov.
\newblock The semigroup of metric measure spaces and its infinitely divisible
  probability measures.
\newblock {\em Trans. Amer. Math. Soc.}, 369(3):1797--1834, 2017.

\bibitem[Eva00]{Ev00}
Steven~N. Evans.
\newblock Kingman's coalescent as a random metric space.
\newblock In {\em Stochastic models ({O}ttawa, {ON}, 1998)}, volume~26 of {\em
  CMS Conf. Proc.}, pages 105--114. Amer. Math. Soc., Providence, RI, 2000.

\bibitem[FG94]{FG94}
Klaus Fleischmann and Andreas Greven.
\newblock Diffusive clustering in an infinite system of hierarchically
  interacting diffusions.
\newblock {\em Probab. Theory Related Fields}, 98(4):517--566, 1994.

\bibitem[FINR04]{FINR04}
Luiz Renato~G. Fontes, Marco Isopi, Charles~M. Newman, and Krishnamurthi
  Ravishankar.
\newblock The {B}rownian web: characterization and convergence.
\newblock {\em Ann. Probab.}, 32(4):2857--2883, 2004.

\bibitem[Fou12]{Fou2012}
Cl{\'e}ment Foucart.
\newblock Generalized {F}leming-{V}iot processes with immigration via
  stochastic flows of partitions.
\newblock {\em ALEA Lat. Am. J. Probab. Math. Stat.}, 9(2):451--472, 2012.

\bibitem[GGR17]{infdiv}
Andreas Greven, Patric~K. Gl\"ode, and Thomas Rippl.
\newblock Branching trees {I}: Concatenation and infinite divisibility.
\newblock {\em ArXive 1612.01265}, submitted June 2017.
\newblock \url{ http://arxiv.org/abs/1612.01265}.

\bibitem[GKW18]{gklimw}
Andreas Greven, Anton Klimovsky, and Anita Winter.
\newblock Evolving genealogies of spatial {$\Lambda$}-{C}annings processes with
  mutation.
\newblock in preparation, 2018.

\bibitem[Gl{\"o}12]{Gl12}
Patric~K. Gl{\"o}de.
\newblock {\em Dynamics of genealogical trees for autocatalytic branching
  processes}.
\newblock PhD thesis, Department Mathematik, Erlangen, Germany, 2012.
\newblock \url{http://nbn-resolving.de/urn:nbn:de:bvb:29-opus-45453}.

\bibitem[GLW05]{GLW05}
Andreas Greven, Vlada Limic, and Anita Winter.
\newblock Representation theorems for interacting {M}oran models, interacting
  {F}isher-{W}right diffusions and applications.
\newblock {\em Electron. J. Probab.}, 10:no. 39, 1286--1356, 2005.

\bibitem[GM18]{gmuk}
Andreas Greven and Chiranjib Mukherjee.
\newblock Genealogies of spatial logistic branching models.
\newblock in preparation, 2018.

\bibitem[GPW09]{GPW09}
Andreas Greven, Peter Pfaffelhuber, and Anita Winter.
\newblock Convergence in distribution of random metric measure spaces
  ({$\Lambda$}-coalescent measure trees).
\newblock {\em Probab. Theory Related Fields}, 145(1-2):285--322, 2009.

\bibitem[GPW13]{GPWmp13}
Andreas Greven, Peter Pfaffelhuber, and Anita Winter.
\newblock Tree-valued resampling dynamics {M}artingale problems and
  applications.
\newblock {\em Probab. Theory Related Fields}, 155(3-4):789--838, 2013.

\bibitem[GRG18]{ggr_GeneralBranching}
Andreas Greven, Thomas Rippl, and Patric~K. Gl\"ode.
\newblock Branching processes — a general concept.
\newblock in preparation, 2018.

\bibitem[Gri17]{Grieshammer2017}
Max Grieshammer.
\newblock {\em {M}easure Representations of Genealogical Processes and
  Applications to {F}leming-{V}iot Models}.
\newblock PhD thesis, Department Mathematik, Erlangen, Germany, 2017.
\newblock \url{http://nbn-resolving.de/urn:nbn:de:bvb:29-opus4-85653}.

\bibitem[GSW]{GSWfoss}
Andreas Greven, Rongfeng Sun, and Anita Winter.
\newblock The evolving genealogy of fossils: Unique characterization by
  martingale problems and applications.
\newblock In preparation 2018.

\bibitem[GSW16]{GSW}
Andreas Greven, Rongfeng Sun, and Anita Winter.
\newblock Continuum space limit of the genealogies of interacting
  {F}leming-{V}iot processes on {$\mathbb Z$}.
\newblock {\em Electron. J. Probab.}, 21:Paper No. 58, 64, 2016.

\bibitem[Guf18]{Gufler2018}
Stephan Gufler.
\newblock Pathwise construction of tree-valued {F}leming-{V}iot processes.
\newblock {\em Electron. J. Probab.}, 23:Paper No. 42, 58, 2018.

\bibitem[KL15]{KL15}
Sandra Kliem and Wolfgang L\"ohr.
\newblock Existence of mark functions in marked metric measure spaces.
\newblock {\em Electron. J. Probab.}, 20:no. 73, 24, 2015.

\bibitem[Kle96]{Kle96}
Achim Klenke.
\newblock Different clustering regimes in systems of hierarchically interacting
  diffusions.
\newblock {\em Ann. Probab.}, 24(2):660--697, 1996.

\bibitem[Kle97]{Kle97}
Achim Klenke.
\newblock Multiple scale analysis of clusters in spatial branching models.
\newblock {\em Ann. Probab.}, 25(4):1670--1711, 1997.

\bibitem[KW18]{KW17}
Sandra Kliem and Anita Winter.
\newblock Evolving phylogenies of trait-dependent branching with mutation and
  competition, part~{I}: Existence.
\newblock {\em Stochastic Process. Appl.}, 2018.

\bibitem[Lab14]{L14}
Cyril Labb\'e.
\newblock From flows of {$\Lambda$}-{F}leming-{V}iot processes to lookdown
  processes via flows of partitions.
\newblock {\em Electron. J. Probab.}, 19:no. 55, 49, 2014.

\bibitem[LG89]{LG89}
Jean-Fran\d{c}ois Le~Gall.
\newblock Marches al\'{e}atoires, mouvement brownien et processus de
  branchement.
\newblock In {\em S\'{e}minaire de {P}robabilit\'{e}s, {XXIII}}, volume 1372 of
  {\em Lecture Notes in Math.}, pages 258--274. Springer, Berlin, 1989.

\bibitem[LG93]{LeGall93}
Jean-Fran\c{c}ois Le~Gall.
\newblock The uniform random tree in a {B}rownian excursion.
\newblock {\em Probab. Theory Related Fields}, 96(3):369--383, 1993.

\bibitem[LG99]{LG99}
Jean-Fran{\c{c}}ois Le~Gall.
\newblock {\em Spatial branching processes, random snakes and partial
  differential equations}.
\newblock Lectures in Mathematics ETH Z\"urich. Birkh\"auser Verlag, Basel,
  1999.

\bibitem[Lig85]{Lig85}
Thomas~M. Liggett.
\newblock {\em Interacting particle systems}, volume 276 of {\em Grundlehren
  der Mathematischen Wissenschaften [Fundamental Principles of Mathematical
  Sciences]}.
\newblock Springer-Verlag, New York, 1985.

\bibitem[L{\"o}h13]{Loehr13}
Wolfgang L{\"o}hr.
\newblock Equivalence of {G}romov-{P}rohorov- and {G}romov's
  {$\underline\square_\lambda$}-metric on the space of metric measure spaces.
\newblock {\em Electron. Commun. Probab.}, 18:no. 17, 10, 2013.

\bibitem[LS81]{LS81}
Thomas~M. Liggett and Frank Spitzer.
\newblock Ergodic theorems for coupled random walks and other systems with
  locally interacting components.
\newblock {\em Z. Wahrsch. Verw. Gebiete}, 56(4):443--468, 1981.

\bibitem[LVW15]{LVW15}
Wolfgang L\"ohr, Guillaume Voisin, and Anita Winter.
\newblock Convergence of bi-measure {$\mathbb{R}$}-trees and the pruning
  process.
\newblock {\em Ann. Inst. Henri Poincar\'e Probab. Stat.}, 51(4):1342--1368,
  2015.

\bibitem[NP89]{NP89}
Jacques Neveu and Jim Pitman.
\newblock The branching process in a {B}rownian excursion.
\newblock In {\em S\'{e}minaire de {P}robabilit\'{e}s, {XXIII}}, volume 1372 of
  {\em Lecture Notes in Math.}, pages 248--257. Springer, Berlin, 1989.

\bibitem[NRS05]{NRS05}
Charles~M. Newman, Krishnamurthi Ravishankar, and Rongfeng Sun.
\newblock Convergence of coalescing nonsimple random walks to the {B}rownian
  web.
\newblock {\em Electron. J. Probab.}, 10:no. 2, 21--60, 2005.

\bibitem[SSS17]{SSS15}
Emmanuel Schertzer, Rongfeng Sun, and Jan~M. Swart.
\newblock The {B}rownian web, the {B}rownian net, and their universality.
\newblock In {\em Advances in disordered systems, random processes and some
  applications}, pages 270--368. Cambridge Univ. Press, Cambridge, 2017.

\bibitem[Win02]{Win02}
Anita Winter.
\newblock Multiple scale analysis of spatial branching processes under the
  {P}alm distribution.
\newblock {\em Electron. J. Probab.}, 7:No. 13, 72, 2002.

\end{thebibliography}

\end{document}